\documentclass{elsart}
\usepackage{amssymb}
\usepackage{amsmath}
\usepackage{graphics}
\usepackage{epsfig}
\usepackage[latin1]{inputenc}

\newtheorem{proposition}{Proposition}

\begin{document}
\begin{frontmatter}
\title{Dual divergence estimators and tests: robustness
results}

\author[label1,label3]{Aida Toma\corauthref{cor1}}, \author[label2]{Michel Broniatowski}

\address[label1]{Mathematics Department, Academy of Economic Studies, Pia\c{t}a Roman\u{a} 6, Bucharest, Romania, e-mail: aida$_{-}$toma@yahoo.com}
\address[label3]{"Gheorghe Mihoc - Caius Iacob" Institute of Mathematical Statistics and Applied Mathematics, Calea 13 Septembrie 13, Bucharest, Romania}
\corauth[cor1]{Corresponding author. Mathematics Department,
Academy of Economic Studies, Pia\c ta Roman\u a 6, Bucharest,
Romania. E-mail: aida$_{-}$toma@yahoo.com}
\address[label2]{LSTA, Université Paris 6, 175 Rue du
Chevaleret, 75013 Paris, France,\\ e-mail:
michel.broniatowski@upmc.fr}

\begin{abstract}
The class of dual $\phi$-divergence estimators (introduced in
Broniatowski and Keziou (2009) \cite{BrKe07}) is explored with
respect to robustness through the influence function approach. For
scale and location models, this class is investigated in terms of
robustness and asymptotic relative efficiency. Some hypothesis
tests based on dual divergence criterions are proposed and their
robustness properties are studied. The empirical performances of
these estimators and tests are illustrated by Monte Carlo
simulation for both noncontaminated and contaminated data.

\end{abstract}


\end{frontmatter}

\section{Introduction}

Minimum divergence estimators and related methods have received
considerable attention in statistical inference because of their
ability to reconcile efficiency and robustness. Among others,
Beran \cite{Be}, Tamura and Boos \cite{TaBo}, Simpson
\cite{Si87,Si89} and Toma \cite{To} proposed families of
parametric estimators minimizing the Hellinger distance between a
nonparametric estimator of the observations density and the model.
They showed that those estimators are both asymptotically
efficient and robust. Generalizing earlier work based on the
Hellinger distance, Lindsay \cite{Li}, Basu and Lindsay
\cite{BaLi}, Morales et al. \cite{MoPaVa} have investigated
minimum divergence estimators, for both discrete and continuous
models. Some families of estimators based on approximate
divergence criterions have also been considered; see Basu et al.
\cite{BaHa}.

Broniatowski and Keziou \cite{BrKe07} have introduced a new
minimum divergence estimation method based on a dual
representation of the divergence between probability measures.
Their estimators are defined in an unified way for both continuous
and discrete models. They do not require any prior smoothing and
include the classical maximum likelihood estimators as a
benchmark. A special case for the Kullback-Leibler divergence is
presented in Broniatowski \cite{Br}. The present paper presents
robustness studies for the classes of estimators generated by the
minimum dual $\phi$-divergence method, as well as for some tests
based on corresponding estimators of the divergence criterion.

We give general results that allow to identify robust estimators
in the class of dual $\phi$-divergence estimators. We apply this
study for the Cressie-Read divergences and state explicit
robustness results for scale models and location models. Gain in
robustness is often paid by some loss in efficiency. This is
discussed for some scale and location models. Our main remarks are
as follows. All the relevant information pertaining to the model
and the true value of the parameter to be estimated should be used
in order to define, when possible, robust and nearly efficient
procedures. Some models allow for such procedures. The example
provided by the scale normal model shows that the choice of a good
estimation criterion is heavily dependent on the acceptable loss
in efficiency in order to achieve a compromise with the robustness
requirement. When sampling under the model is overspread
(typically for Cauchy and logistic models), non surprisingly the
maximum likelihood estimator is both efficient and robust and
therefore should be prefered (see subsection \ref{study}).

On the other hand, these estimation results constitute the
premises to construct some robust tests. The purpose of robust
testing is twofold. First, the level of a test should be stable
under small arbitrary departures from the null hypothesis (i.e.
robustness of validity). Second, the test should have a good power
under small arbitrary departures from specified alternatives (i.e.
robustness of efficiency). To control the test stability against
outliers in the aforementioned senses, we compute the asymptotic
level of the test under a sequence of contaminated null
distributions, as well as the asymptotic power of the test under a
sequence of contaminated alternatives. These quantities are seen
to be controlled by the influence function of the test statistic.
In this way, the robustness of the test is a consequence of the
robustness of the test statistic based on a dual $\phi$-divergence
estimator. In many cases, this requirement is met when the dual
$\phi$-divergence estimator itself is robust.

The paper is organized as follows: in Section 2 we present the
classes of estimators generated by the minimum dual
$\phi$-divergence method. In Section 3, for these estimators, we
compute the influence functions and give the Fisher consistency.
We particularize this study for the Cressie-Read divergences and
state robustness results for scale models and location models.
Section 4 is devoted to hypothesis testing. We give general
convergence results for contaminated observations and use it to
compute the asymptotic level and the asymptotic power for the
tests that we propose. In Section 5, the performances of the
estimators and tests are illustrated by Monte Carlo simulation
studies. In Section \ref{secchoice} we shortly presents a proposal
for the adaptive choice of tuning parameters.

\section{Minimum divergence estimators}
\subsection{Minimum divergence estimators}
Let $\varphi$ be a non-negative convex function defined from
$(0,\infty)$ onto $[0,\infty]$ and satisfying $\varphi(1)=0$. Also
extend $\varphi$ at 0 defining $\varphi(0)=\underset{x\downarrow
0}{\lim} \;\varphi(x)$.  Let $(\mathcal{X}, \mathcal{B})$ be a
measurable space and $P$ be a probability measure (p.m.) defined
on $(\mathcal{X}, \mathcal{B})$. Following Rüschendorf \cite{Ru},
for any p.m. $Q$ absolutely continuous (a.c.) w.r.t. $P$, the
$\phi$-divergence between $Q$ and $P$ is defined by
\begin{equation}\label{eq2.1}
\phi(Q,P):=\int \varphi\left(\frac{dQ}{dP}\right)dP.
\end{equation}
When $Q$ is not a.c. w.r.t. $P,$ we set $\phi(Q,P)=\infty$. We
refer to Liese and Vajda \cite{LiVa06} for an overview on the
origin of the concept of divergence in Statistics.

A commonly used family of divergences is the so-called "power
divergences", introduced by Cressie and Read \cite{CrRe} and
defined by the class of functions
\begin{equation}\label{1.3}
x\in \mathbb{R_{+}^{*}}\mapsto \varphi_{\gamma}(x):=\frac{x^{\gamma}-\gamma x+ \gamma-1}{\gamma(\gamma-1)}
\end{equation}
for $\gamma\in \mathbb{R}\setminus\{ 0,1\}$ and
$\varphi_{0}(x):=-\log x+ x-1$, $\varphi_{1}(x):=x\log x- x+1$
with $\varphi_{\gamma}(0)=\underset{x\downarrow
0}{\lim}\;\varphi_{\gamma}(x)$,
$\varphi_{\gamma}(\infty)=\underset{x\to
\infty}{\lim}\;\varphi_{\gamma}(x)$, for any $\gamma\in
\mathbb{R}.$ The Kullback-Leibler divergence (KL) is associated
with $\varphi_{1}$, the modified Kullback-Leibler (KL$_{m}$) to
$\varphi_{0}$, the $\chi^{2}$ divergence to $\varphi_{2}$, the
modified $\chi^{2}$ divergence ($\chi_{m}^{2}$) to $\varphi_{-1}$
and the Hellinger distance to $\varphi_{1/2}$.

Let $\{P_{\theta}: \theta\in \Theta\}$ be some identifiable
parametric model with $\Theta$ a subset of $\mathbb{R}^{d}$.
Consider the problem of estimation of the unknown true value of
the parameter $\theta_{0}$ on the basis of an i.i.d. sample
$X_{1}, \dots, X_{n}$ with p.m. $P_{\theta_{0}}$.

When all p.m. $P_{\theta}$ share the same finite support $S$ which is independent upon the parameter $\theta$, the
$\phi$-divergence between $P_{\theta}$ and $P_{\theta_{0}}$ has the form
\begin{equation*}
\phi(P_{\theta},P_{\theta_{0}})=\sum_{j\in
S}\varphi\left(\frac{P_{\theta}(j)}{P_{\theta_{0}}(j)}\right)P_{\theta_{0}}(j).
\end{equation*}
In this case, Liese and Vajda \cite{LiVa87}, Lindsay \cite{Li} and
Morales et al. \cite{MoPaVa} investigated the so-called "minimum
$\phi$-divergence estimators" (minimum disparity estimators in
Lindsay \cite{Li}) of the parameter $\theta_{0}$ defined by
\begin{equation}\label{eq1}
\widehat{\theta}_{n}:=\arg\inf_{\theta\in
\Theta}\phi(P_{\theta},P_{n}),
\end{equation}
where $\phi(P_{\theta},P_{n})$ is the plug-in estimator of $\phi(P_{\theta},P_{\theta_{0}})$
\begin{equation*}
\phi(P_{\theta},P_{n})=\sum_{j\in
S}\varphi\left(\frac{P_{\theta}(j)}{P_{n}(j)}\right)P_{n}(j),
\end{equation*}
$P_{n}$ being the empirical measure associated to the sample. The
interest on these estimators is motivated by the fact that a
suitable choice of the divergence may leads to an estimator more
robust than the maximum likelihood one (see also Jiménez and Shao
\cite{JiSh}). For continuous models, the estimators in (\ref{eq1})
are not defined. Basu and Lindsay \cite{BaLi}, among others,
proposed smoothed versions of (\ref{eq1}) in this case.

In the following, for notational clearness we write
$\phi(\alpha,\theta)$ for $\phi(P_{\alpha},P_{\theta})$ for
$\alpha$ and $\theta$ in $\Theta$. We assume that for any
$\theta\in \Theta$, $P_{\theta}$ has density $p_{\theta}$ with
respect to some dominating $\sigma$-finite measure $\lambda$.

The divergence $\phi(\alpha,\theta_0)$ can be represented as
resulting from an optimization procedure. This result has been
obtained independently by Liese and Vajda \cite{LiVa06} and
Broniatowski and Keziou \cite{BrKe06} who called it the dual form
of a divergence, due to its connection with convex analysis.
Assuming the strict convexity and the differentiability of the
function $\varphi$, it holds
\begin{equation}\label{incon}
\varphi(t)\geq\varphi(s)+\varphi^{\prime}(s)(t-s)
\end{equation}
with equality only for $s=t$. Let $\alpha$ and $\theta_0$ be fixed
and put $t=p_{\alpha}(x)/p_{\theta_0}(x)$ and
$s=p_{\alpha}(x)/p_{\theta}(x)$ in (\ref{incon}) and then
integrate with respect to $P_{\theta_0}$. This gives
\begin{equation}\label{drep}
\phi(\alpha, \theta_{0})=\int
\varphi\left(\frac{p_{\alpha}}{p_{\theta_0}}\right)dP_{\theta_0}=\sup_{\theta\in
\Theta}\int m(\theta,\alpha)dP_{\theta_{0}}
\end{equation}
with $m(\theta,\alpha):x\mapsto m(\theta,\alpha,x)$ and
\begin{equation}\label{e2b}
m(\theta,\alpha,x):=\int
\varphi^{\prime}\left(\frac{p_{\alpha}}{p_{\theta}}\right)dP_{\alpha}-\left\{\varphi^{\prime}\left(\frac{p_{\alpha}}{p_{\theta}}(x)\right)\frac{p_{\alpha}}{p_{\theta}}(x)-\varphi\left(\frac{p_{\alpha}}{p_{\theta}}(x)\right)\right\}.
\end{equation}

The supremum in (\ref{drep}) is unique and is attained in
$\theta=\theta_0$, independently upon the value of $\alpha$.
Naturally, a class of estimators of $\theta_{0}$, called "dual
$\phi$-divergence estimators" (D$\phi$E's), is defined by
\begin{equation}\label{f2.15}
\widehat{\theta}_{n}(\alpha):=\arg\sup_{\theta\in \Theta}\int
m(\theta,\alpha)dP_{n},\;\;\alpha\in \Theta.
\end{equation}
Formula (\ref{f2.15}) defines a family of M-estimators indexed by
some instrumental value of the parameter $\alpha$ and by the
function $\varphi$ defining the divergence. The choice of $\alpha$
appears as a major feature in the estimation procedure. Its value
is strongly dependent upon some a priori knowledge on the value of
the parameter to be estimated. In some examples in subsection
\ref{study}, it even appears that a sharp a priori knowledge on
the order of $\theta_{0}$ leads to nearly efficient and robust
estimates. This plays in favor of using the available information
pertaining to the model and the data. Section \ref{secchoice}
shortly presents some proposal for the adaptive choice of
$\alpha$.

For each $\alpha\in \Theta$, the divergence
$\phi(P_{\alpha},P_{\theta_{0}})$ between $P_{\alpha}$ and
$P_{\theta_{0}}$ is estimated by
\begin{equation}\label{ff2.15}
\widehat{\phi}_{n}(\alpha,\theta_{0}):=\int
m(\widehat{\theta}_{n}(\alpha),\alpha)dP_{n}=\sup_{\theta\in
\Theta}\int m(\theta,\alpha)dP_{n}.
\end{equation}

Further, since
\begin{equation*}
\inf_{\alpha\in\Theta}\phi(\alpha,\theta_{0})=\phi(\theta_{0},\theta_{0})=0,
\end{equation*}
and since the infimum in the above display is unique due to the
strict convexity of $\varphi$, a natural definition of estimators
of $\theta_{0}$, called "minimum dual $\phi$-divergence
estimators" (MD$\phi$E's), is provided by
\begin{equation}\label{fff2.15}
\widehat{\alpha}_{n}:=\arg\inf_{\alpha\in\Theta}\widehat{\phi}_{n}(\alpha,\theta_{0})=\arg\inf_{\alpha\in\Theta}\sup_{\theta\in
\Theta}\int m(\theta,\alpha)dP_{n}.
\end{equation}
The D$\phi$E's enjoy the same invariance property as the maximum
likelihood estimator does. Invariance with respect to a
reparametrization (one to one transformation of the parameter
space) holds with direct substitution in (\ref{f2.15}). Also,
consider a one to one differentiable transformation of the
observations, say $Y=T(X)$ and the Jacobian
$J(x)=\frac{d}{dx}T(x)$. Let $\widehat{\theta}_{n}(\alpha)$
defined in (\ref{f2.15}), based on the $X_{i}$'s. Let
$f_{\theta}(y)$ denote the density of the transformed variable $Y$
and $\widehat{\theta}^{*}_{n}(\alpha)$ be the D$\phi$E based on
the $Y_{i}$'s in the transformed model (with the same parameter
$\theta$). Specifically,
\begin{equation*}
\widehat{\theta}^{*}_{n}(\alpha)=\arg\sup_{\theta\in\Theta}\left\{\int\varphi^{'}\left(\frac{f_{\alpha}}{f_{\theta}}(y)\right)f_{\alpha}(y)dy-\frac{1}{n}\sum_{i=1}^{n}\left(\varphi^{'}\left(\frac{f_{\alpha}}{f_{\theta}}(Y_{i})\right)\frac{f_{\alpha}}{f_{\theta}}(Y_{i})-\varphi\left(\frac{f_{\alpha}}{f_{\theta}}(Y_{i})\right)\right)\right\}.
\end{equation*}
Since
\begin{equation*}
f_{\theta}(y)=p_{\theta}(T^{-1}(y))|J(T^{-1}(y))|^{-1}
\end{equation*}
for all $\theta\in\Theta$, it follows that
$\widehat{\theta}^{*}_{n}(\alpha)=\widehat{\theta}_{n}(\alpha)$,
which is to say that the D$\phi$E's are invariant estimators under
any regular transformation of the observation space. The same
invariance properties hold for MD$\phi$E's.

Broniatowski and Keziou \cite{BrKe07} have proved both the weak
and the strong consistency, as well as the asymptotic normality
for the estimators $\widehat{\theta}_{n}(\alpha)$ and
$\widehat{\alpha}_{n}$. In the next sections, we study robustness
properties for these classes of estimators and robustness of some
tests based on dual $\phi$-divergence estimators.
\subsection{Some comments on robustness}
The special form of divergence based estimators to be studied in
this paper leads us to handle robustness characteristics through
the influence function approach. An alternative and appealing
robustness analysis in the minimum divergence methods is provided
by the Residual Adjustment Function (RAF) (introduced in Lindsay
\cite{Li}), which explains the incidence of non typical Pearson
residuals, corresponding to over or sub-sampling, in the stability
of the estimates. This method is quite natural for finitely
supported models. In the case when the densities in the model are
continuous, the Pearson residuals are estimated non parametrically
which appears to cause quite a number of difficulties when adapted
to minimum dual divergence estimation. This motivates the present
choice in favor of the influence function approach.

Let $\alpha$ be fixed. For the Cressie-Read divergences, the
equation whose solution is $\widehat{\theta}_{n}(\alpha)$ defined
by (\ref{f2.15}) is
\begin{equation}\label{equa}
-\int
\left(\frac{p_{\alpha}}{p_{\theta}}\right)^{\gamma}\dot{p}_{\theta}d\lambda+\frac{1}{n}\sum_{i=1}^{n}\left(\frac{p_{\alpha}(X_{i})}{p_{\theta}(X_{i})}\right)^{\gamma}\frac{\dot{p}_{\theta}(X_{i})}{p_{\theta}(X_{i})}=0,
\end{equation}
where $\dot{p}_{\theta}$ is the derivative with respect to
$\theta$ of $p_{\theta}$. Starting from the definition given by
(\ref{f2.15}), this equation is obtained by equalizing with zero
the derivative with respect to $\theta$ of $\int
m(\theta,\alpha)dP_{n}$.

Let $x$ be some outlier. The role of $x$ in (\ref{equa}) is handled in the term
\begin{equation}\label{term}\left(\frac{p_{\alpha}(x)}{p_{\theta}(x)}\right)^{\gamma}\frac{\dot{p}_{\theta}(x)}{p_{\theta}(x)}.\end{equation}
 The more stable this term, the more robust the estimate. In the classical case of the maximum likelihood estimator
 (which corresponds to $\widehat{\theta}_{n}(\alpha)$ with $\gamma=0$ and independent on $\alpha$),
 this term writes as $\frac{\dot{p}_{\theta}(x)}{p_{\theta}(x)}$ which is the likelihood score function associated to $x$.
 It is well known that, for most models, this term is usually unbounded when $x$ belongs to $\mathbb{R}$, saying that
  the maximum likelihood estimator is not robust. In this respect, (\ref{term}) appears as a weighted likelihood score function.
  In our approach, for several models, such as the normal scale, (\ref{term}) is a bounded function of $x$, although
  $\frac{\dot{p}_{\theta}(x)}{p_{\theta}(x)}$ itself is not. Thus, in estimating equation (\ref{equa}), the score function is downweighted for
  large observations. The robustness of $\widehat{\theta}_{n}(\alpha)$ comes as a downweight effect of the
  quantity $\frac{\dot{p}_{\theta}(X_{i})}{p_{\theta}(X_{i})}$ through the multiplicative term $\left(\frac{p_{\alpha}(X_{i})}{p_{\theta}(X_{i})}\right)^{\gamma}$
  which depends on the choice of the divergence. This choice is dictated by the form of $\frac{p_{\alpha}(x)}{p_{\theta}(x)}$ for large $x$ and $\alpha$ fixed.
  For the models we'll consider as examples, for large $x$ and $\alpha$ fixed, the quantity $\frac{p_{\alpha}(x)}{p_{\theta}(x)}$ can be large,
  close to zero, or close to one. Then we appropriately choose $\gamma$ to be negative, respectively positive in order to obtain the downweight effect.
  In the next section we study in detail these robustness properties by the means of the influence function.

Some alternative choice has been proposed in literature. Basu et
al. \cite{BaHa} proposed to alter the likelihood score factor by
the multiplicative term $p_{\theta}^{\beta}(x)$, where $\beta>0$.
This induces an estimating procedure which is connected to the
minimization of a density power divergence. Both their approach
and the present one are adaptive in the sense that the downweight
likelihood score factor is calibrated on the data.

Robustness as handled in the present paper is against the bias due
to the presence of very few outliers in the data set. Bias due to
misspecification of the model is not considered. It has been
observed that D$\phi$E's are biased under misspecification even in
simple situations (for example when estimating the mean in a
normal model with assumed variance 1, whereas the true variance is
not 1); see Broniatowski and Vajda \cite{BrVa09}; similar bias are
unavoidable in parametric inference and can only be reduced
through adaptive specific procedures, not studied here. For
alternative robust M-estimation methods using divergences we refer
to Toma \cite{To3}.

\section{Robustness of the estimators}
\subsection{Fisher consistency and influence functions}
In order to measure the robustness of an estimator it is common to compute the influence function of the
corresponding functional.

A map $T$ which
sends an arbitrary probability measure into the parameter space is a statistical functional corresponding to an
estimator $T_{n}$ of the parameter $\theta$ whenever $T(P_{n})=T_{n}$.

This functional is called Fisher consistent for the parametric model $\{P_{\theta}:\theta\in \Theta\}$ if
$T(P_{\theta})=\theta$, for all $\theta\in\Theta$.

The influence function of the functional $T$ in $P$ measures the effect on $T$ of adding a small
mass at $x$ and is defined as
\begin{equation}
\label{e1b} \mathrm{IF}(x;T,P)=\lim_{\varepsilon \rightarrow 0}\frac{T(\widetilde{P}_{\varepsilon
x})-T(P)}{\varepsilon}
\end{equation}
where $\widetilde{P}_{\varepsilon
x}=(1-\varepsilon)P+\varepsilon\delta_{x}$ and $\delta_{x}$ is the
Dirac measure putting all its mass at $x$.

The gross error sensitivity measures approximately the maximum
contribution to the estimation error that can be produced by a
single outlier and is defined as
\begin{equation*}\sup_{x}\|\mathrm{IF}(x;T,P)\|.\end{equation*}
Whenever the gross error sensitivity is finite, the estimator
associated with the functional $T$ is called B-robust.

Let
$X_{1},\dots, X_{n}$ be an i.i.d. sample with p.m. $P$.

Let $\alpha$ be fixed and consider the dual $\phi$-divergence
estimators $\widehat{\theta}_{n}(\alpha)$ defined in (\ref{f2.15}). The functional associated to an estimator $\widehat{\theta}_{n}(\alpha)$ is
\begin{equation}
\label{e3b} T_{\alpha}(P):= \arg\sup_{\theta\in\Theta}\int m(\theta,\alpha,y)dP(y).
\end{equation}

The functional $T_{\alpha}$ is Fisher consistent. Indeed, the
function $\theta \mapsto \int m(\theta,\alpha)dP_{\theta_{0}}$ has
a unique maximizer $\theta=\theta_{0}$. Therefore
$T_{\alpha}(P_{\theta})=\theta,$ for all $\theta\in\Theta$.

We denote $m^{\prime}(\theta,\alpha)=\frac{\partial}{\partial
\theta}m(\theta,\alpha)$ the $d$-dimensional column vector with
entries $\frac{\partial}{\partial \theta_{i}}m(\theta,\alpha)$ and
$m^{\prime\prime}(\theta,\alpha)$ the $d\times d$ matrix with
entries
$\frac{\partial^{2}}{\partial\theta_{i}\partial\theta_{j}}m(\theta,\alpha)$.

In the rest of the paper, for each $\alpha$, we suppose that the
function $\theta\mapsto m(\theta,\alpha)$ is twice continuously
differentiable and that the matrix $\int
m^{\prime\prime}(\theta_{0},\alpha)dP_{\theta_{0}}$ exists and is
invertible. We also suppose that, for each $\alpha$, all the
partial derivatives of order 1 and 2 of the function
$\theta\mapsto m(\theta,\alpha)$ are respectively dominated on
some neighborhoods of $\theta_{0}$ by $P_{\theta_{0}}$-integrable
functions. This justifies the subsequent interchanges of
derivation with respect to $\theta$ and integration.
\begin{proposition}
\label{t1b} The influence function of the functional $T_{\alpha}$ corresponding to an estimator
$\widehat{\theta}_{n}(\alpha)$ is given by
\begin{eqnarray*}
\label{e4b} \mathrm{IF}(x;T_{\alpha},P_{\theta_{0}}) &=& \left[\int
m^{\prime\prime}(\theta_{0},\alpha)dP_{\theta_{0}}\right]^{-1}\left\{\int
\varphi^{\prime\prime}\left(\frac{p_{\alpha}}{p_{\theta_{0}}}\right)\frac{p_{\alpha}}{p_{\theta_{0}}^{2}}\dot{p}_{\theta_{0}}dP_{\alpha}-\right.\\
&&
-\left.\varphi^{\prime\prime}\left(\frac{p_{\alpha}}{p_{\theta_{0}}}(x)\right)\frac{p_{\alpha}^{2}(x)}{p_{\theta_{0}}^{3}(x)}\dot{p}_{\theta_{0}}(x)\right\}.
\end{eqnarray*}
\end{proposition}

Particularizing $\alpha=\theta_{0}$ in Proposition \ref{t1b} yields
\begin{equation*}
\mathrm{IF}(x;T_{\theta_{0}},P_{\theta_{0}}) = -\left[\int
m^{\prime
\prime}(\theta_{0},\theta_{0})dP_{\theta_{0}}\right]^{-1}\varphi^{\prime
\prime}(1)\frac{\dot{p}_{\theta_{0}}(x)}{p_{\theta_{0}}(x)}\end{equation*}
and taking into account that
\begin{equation*}
-\left[\int m^{\prime
\prime}(\theta_{0},\theta_{0})dP_{\theta_{0}}\right]^{-1}=\frac{1}{\varphi^{\prime\prime}(1)}I_{\theta_{0}}^{-1}
\end{equation*}
it holds
\begin{equation}
\label{e6b}
\mathrm{IF}(x;T_{\theta_{0}},P_{\theta_{0}})=I_{\theta_{0}}^{-1}\frac{\dot{p}_{\theta_{0}}(x)}{p_{\theta_{0}}(x)}
\end{equation}
where $I_{\theta_{0}}$ is the information matrix
$I_{\theta_{0}}=\int\frac{\dot{p}_{\theta_{0}}\dot{p}_{\theta_{0}}^{t}}{p_{\theta_{0}}}d \lambda$.

We now look at the corresponding estimators of the
$\phi$-divergence. For fixed $\alpha$, the divergence
$\phi(P_{\alpha},P)$ between the probability measures $P_{\alpha}$
and $P$ is estimated by (\ref{ff2.15}). The statistical functional
associated to $\widehat{\phi}_{n}(P_{\alpha},P_{\theta_{0}})$ is
\begin{equation}
\label{e7b} U_{\alpha}(P):=\int m(T_{\alpha}(P),\alpha,y)dP(y).
\end{equation}

The functional $U_{\alpha}$ has the property that
$U_{\alpha}(P_{\theta})=\phi(\alpha,\theta)$, for any
$\theta\in\Theta$. Indeed, using the fact that $T_{\alpha}$ is a
Fisher consistent functional,
\begin{equation*}
U_{\alpha}(P_{\theta})= \int m\left(T_{\alpha}(P_{\theta}\right),\alpha,y)dP_{\theta}(y)=\int
m(\theta,\alpha,y)dP_{\theta}(y)= \phi(\alpha,\theta)
\end{equation*}
for all $\theta\in\Theta$.

\begin{proposition}
\label{t2b} The influence function of the functional $U_{\alpha}$
corresponding to the estimator $\widehat{\phi}_{n}(P_{\alpha},P)$
is given by
\begin{equation}
\label{e8b} \mathrm{IF}(x;U_{\alpha},P_{\theta_{0}})=-\phi(\alpha,\theta_{0})+m(\theta_{0},\alpha,x).
\end{equation}
\end{proposition}

For a minimum dual $\phi$-divergence estimator $\widehat{\alpha}_{n}$ defined in (\ref{fff2.15}), the corresponding functional is
\begin{equation}
\label{e10b} V(P):=\arg\inf_{\alpha\in\Theta}U_{\alpha}(P)=\arg\inf_{\alpha\in\Theta}\int
m(T_{\alpha}(P),\alpha,y)dP(y).
\end{equation}

The statistical functional $V$ is Fisher consistent. Indeed,
\begin{equation*}
V(P_{\theta})=\arg\inf_{\alpha\in\Theta}U_{\alpha}(P_{\theta})=\arg\inf_{\alpha\in\Theta}\phi(\alpha,\theta)=\theta
\end{equation*}
for all $\theta\in \Theta$.

In the following proposition, we suppose that the function
$m(\theta,\alpha)$ admits partial derivatives of order 1 and 2
with respect to $\theta$ and $\alpha$ and also we suppose that conditions permitting to derivate $m(\theta,\alpha)$
under the integral sign hold. The following result states that, unlike
$\widehat{\theta}_{n}(\alpha)$, an estimator
$\widehat{\alpha}_{n}$ is generally not robust. Indeed, it has the
same robustness properties as the maximum likelihood estimator,
since it has its influence function which in most cases is unbounded. Whatever the divergence, the estimators $\widehat{\alpha}_{n}$ have
the same influence function.
\begin{proposition}\label{t3b} The influence function of the functional $V$ corresponding to an
estimator $\widehat{\alpha}_{n}$ is given by
\begin{equation}
\label{e11b} \mathrm{IF}(x;V,P_{\theta_{0}})=I_{\theta_{0}}^{-1}\frac{\dot{p}_{\theta_{0}}(x)}{p_{\theta_{0}}(x)}.
\end{equation}
\end{proposition}

\subsection{Robustness of the estimators for scale models and location models}\label{study}In this subsection, examining the
expressions of the influence functions, we give conditions for
attaining the B-robustness of the dual $\phi$-divergence
estimators $\widehat{\theta}_{n}(\alpha)$, as well as of the
corresponding divergence estimators. The case of interest in our
B-robustness study is $\alpha\neq\theta_{0}$ since, as observed
above, the choice $\alpha=\theta_{0}$ generally leads to unbounded
influence functions. For the Cressie-Read family of divergences
(\ref{1.3}) it holds
\begin{equation}
\label{b1}\mathrm{IF}(x;T_{\alpha}, P_{\theta_{0}})=[\int
m^{\prime\prime}(\theta_{0},\alpha)dP_{\theta_{0}}]^{-1}\left\{\int\left(\frac{p_{\alpha}}{p_{\theta_{0}}}\right)^{\gamma}\dot{p}_{\theta_{0}}d\lambda-\left(\frac{p_{\alpha}(x)}{p_{\theta_{0}}(x)}\right)^{\gamma}\frac{\dot{p}_{\theta_{0}}(x)}{p_{\theta_{0}}(x)}\right\}
\end{equation}
and
\begin{eqnarray*}
\mathrm{IF}(x;U_{\alpha},P_{\theta_{0}})&=&-\phi
(\alpha,\theta_{0})+m(\theta_{0},\alpha,x)\\ &=& -\phi
(\alpha,\theta_{0})+\frac{1}{\gamma-1}\left\{\int
\left(\frac{p_{\alpha}}{p_{\theta_{0}}}\right)^{\gamma-1}\!\!dP_{\alpha}-1\right\}-\frac{1}{\gamma}\left\{
\left(\frac{p_{\alpha}(x)}{p_{\theta_{0}}(x)}\right)^{\gamma}\!\!-1\right\}.
\end{eqnarray*}

\subsubsection{Scale models}\label{scalem}
For a given density $p$, it holds $p_{\theta}(x)=\frac{1}{\theta}p(\frac{x}{\theta})$ and
$\dot{p}_{\theta}(x)=-\frac{1}{\theta^{2}}\left[p\left(\frac{x}{\theta}\right)+\frac{x}{\theta}\dot{p}\left(\frac{x}{\theta}\right)\right]$. Consider the following conditions:

(A.1) $\int |u  \dot{p}(u)|du< \infty$.

(A.2)
$\sup_{x}\frac{p(\alpha^{-1}x)}{p(\theta_{0}^{-1}x)}<\infty$.

(A.3)
$\sup_{x}\frac{p(\theta_{0}^{-1}x)}{p(\alpha^{-1}x)}<\infty$.

(A.4) $\sup_{x}\left|\frac{\partial}{\partial\theta}[\log
p(\theta_{0}^{-1}x)]\left(\frac{p(\alpha^{-1}x)}{p(\theta_{0}^{-1}x)}\right)^{\gamma}\right|<\infty$.


\begin{proposition}\label{psc}
For scale models, if the conditions $\mathrm{(A.2)}$
$\mathrm{(}$for the case $\gamma>0\mathrm{)}$ or $\mathrm{(A.3)}$
$\mathrm{(}$for the case $\gamma<0\mathrm{)}$ together with
$\mathrm{(A.1)}$ and $\mathrm{(A.4)}$ are satisfied, then
$\widehat{\theta}_{n}(\alpha)$ is B-robust.
\end{proposition}

As a particular case, consider the problem of robust estimation of
the parameter $\theta_{0}=\sigma$ of the univariate normal model,
when the mean $m$ is known, intending to use an estimator
$\widehat{\theta}_{n}(\overline{\sigma})$ with
$\overline{\sigma}\neq \sigma$. We are interested on those
divergences from the Cressie-Read family and those possible values
of $\overline{\sigma}$ for which
$\widehat{\theta}_{n}(\overline{\sigma})$ is B-robust. We have
\begin{equation*}
\mathrm{IF}(x;T_{\overline{\sigma}},P_{\sigma})=[\int
m^{\prime\prime}(\sigma,\overline{\sigma})dP_{\sigma}]^{-1}\left\{\int\left(\frac{p_{\overline{\sigma}}}{p_{\sigma}}\right)^{\gamma}\frac{\dot{p}_{\sigma}}{p_{\sigma}}dP_{\sigma}-\left(\frac{p_{\overline{\sigma}}(x)}{p_{\sigma}(x)}\right)^{\gamma}\frac{\dot{p}_{\sigma}(x)}{p_{\sigma}(x)}\right\}.
\end{equation*}

\begin{figure}[h!]
\begin{center}
\includegraphics[width=8cm,height=8cm]{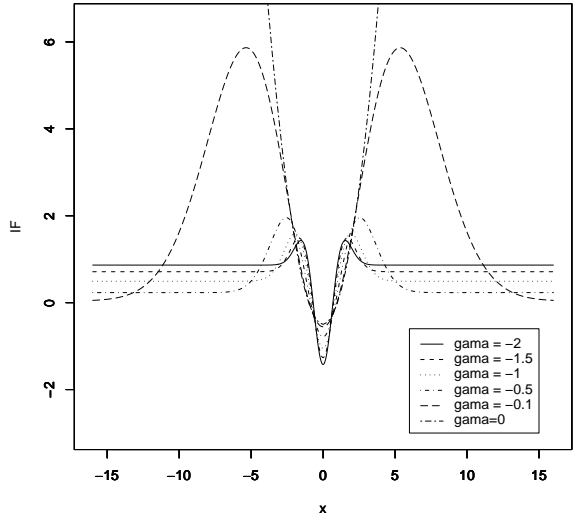}
\end{center}
\caption{Influence functions
$\mathrm{IF}(x;T_{\overline{\sigma}},P_{\sigma})$ for normal scale
model, when $m=0$, the true scale parameter is $\sigma=1$ and
$\overline{\sigma}=1.9$. }\label{BOIF}
\end{figure}

It is easily seen that
$\mathrm{IF}(x;T_{\overline{\sigma}},P_{\sigma})$ is bounded
whenever the function
$\left(\frac{p_{\overline{\sigma}}(x)}{p_{\sigma}(x)}\right)^{\gamma}\frac{\dot{p}_{\sigma}(x)}{p_{\sigma}(x)}$
is bounded. Since
\begin{equation}\label{eqd1}
\left(\frac{p_{\overline{\sigma}}(x)}{p_{\sigma}(x)}\right)^{\gamma}\frac{\dot{p}_{\sigma}(x)}{p_{\sigma}(x)}=\frac{\sigma^{\gamma-1}}{\overline{\sigma}^{\gamma}}\left\{\left(\frac{x-m}{\sigma}\right)^{2}-1\right\}\left(\exp\left(-\frac{1}{2}\left\{\left(\frac{x-m}{\overline{\sigma}}\right)^{2}-\left(\frac{x-m}{\sigma}\right)^{2}\right\}\right)\right)^{\gamma}
\end{equation}
boundedness of $\mathrm{IF}(x;T_{\overline{\sigma}},P_{\sigma})$
holds when $\gamma>0$ and $\overline{\sigma}<\sigma$ or when
$\gamma<0$ and $\overline{\sigma}>\sigma$, cases in which the
conditions of Proposition \ref{psc} are satisfied. A simple
calculation shows that these choices of $\gamma$ and
$\overline{\sigma}$ assure that $\int
m^{\prime\prime}(\sigma,\overline{\sigma})dP_{\sigma}$ is finite
and non zero. However, when using the modified Kullback-Leibler
divergence ($\gamma$=0), none of the estimators
$\widehat{\theta}_{n}(\overline{\sigma})$ is B-robust, the
function (\ref{eqd1}) being unbounded. These aspects can also be
observed in Figure \ref{BOIF}, which presents influence functions
for different divergences when $\sigma=1$ and
$\overline{\sigma}=1.9$. The negative values of the influence
function in a neighborhood of 0 is explained by the decrease of
the variance estimate when oversampling close to the mean.

The asymptotic relative efficiency of an estimator is the ratio of
the asymptotic variance of the maximum likelihood estimator to
that of the estimator in question. For the scale normal model, the
choice of $\overline{\sigma}$ close to $\sigma$ assures a good
efficiency of $\widehat{\theta}_{n}(\overline{\sigma})$ and also
the B-robustness property. Then, the bigger is the value of
$|\gamma|$, the smaller is the gross error sensitivity of the
estimator. For example, for $\sigma=1$ and
$\overline{\sigma}=0.99$, the efficiency of
$\widehat{\theta}_{n}(\overline{\sigma})$ is 0.9803 when
$\gamma=0.5$, 0.9615 when $\gamma=1$, 0.9266 when $\gamma=2$ and
0.8947 when $\gamma=3$, the most B-robust estimator corresponding
to $\gamma=3$. As can be inferred from Figure \ref{BOIF}, the
curves $\mathrm{IF}^2(x;T_{\overline{\sigma}},P_{\sigma})$ are
ordered decreasingly with respect to $|\gamma|$. Therefore, large
values of $|\gamma|$ lead to small gross error sensitivities and
low efficiencies, since the asymptotic variance of
$\widehat{\theta}_{n}(\overline{\sigma})$ is $[\int
\mathrm{IF}^2(x;T_{\overline{\sigma}},P_{\sigma})dP_{\sigma}]^{-1}$
(see also Hampel et al. \cite{HaRoRoSt} for this formula).

For scale models, conditions of Proposition \ref{psc} assure that
$\widehat{\theta}_{n}(\alpha)$ and the corresponding divergence
estimator $\widehat{\phi}_{n}(\alpha,\theta_{0})$ are B-robust.

\subsubsection{Location models}\label{sublocmo}

It holds $p_{\theta}(x)=p(x-\theta)$.
\begin{proposition}\label{loccondd}
For location models, if the condition
\begin{equation}\label{loccond}
\sup_{x}\left|\left(\frac{p(x-\alpha)}{p(x-\theta_{0})}\right)^{\gamma}\frac{\partial}{\partial\theta}\log
p(x-\theta_{0})\right|<\infty
\end{equation}
is satisfied, then $\widehat{\theta}_{n}(\alpha)$ is B-robust.
\end{proposition}

\begin{figure}[]
\begin{center}
\includegraphics[width=8cm,height=8cm]{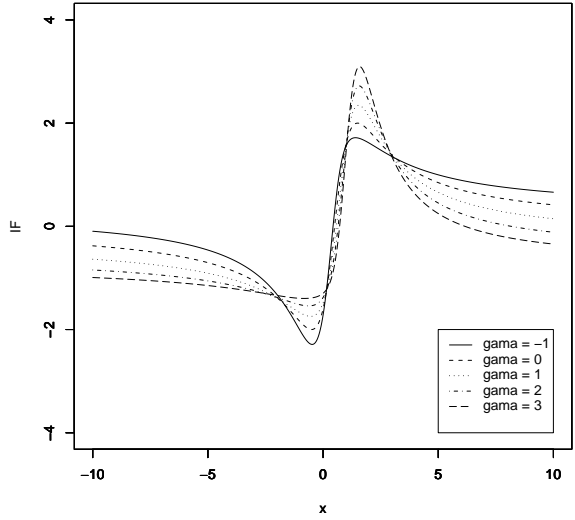}
\end{center}
\caption{Influence functions
$\mathrm{IF}(x;T_{\alpha},P_{\theta_{0}})$ for the Cauchy location
model, when the true location parameter is $\theta_{0}=0.5$ and
$\alpha=0.8$. }\label{Co1}
\end{figure}

For the Cauchy density the maximum likelihood estimator exists,
it is consistent, efficient and B-robust and all the estimators
$\widehat{\theta}_{n}(\alpha)$ exist and are B-robust. Indeed, condition (\ref{loccond}) writes
\begin{equation*}
\sup_{x}2\left|\left(\frac{1+(x-\theta_{0})^{2}}{1+(x-\alpha)^{2}}\right)^{\gamma}\frac{x-\theta_{0}}{1+(x-\theta_{0})^{2}}\right|<\infty
\end{equation*}
and is fulfilled for any $\gamma$ and any $\alpha$. Also, the
integral $\int m^{\prime\prime}(\theta_{0},\alpha)dP_{\theta_{0}}$
exists and is different to zero for any $\gamma$ and any $\alpha$.
This is quite natural since sampling of the Cauchy law makes
equivalent outliers and large sample points due to heavy tails.
However it is known that the likelihood equation for Cauchy
distribution has multiple roots. The number of solutions behaves
asymptotically as two times a Poisson(1/$\pi$) variable plus 1
(see van der Vaart \cite{va} p. 74). The possible selection rule
for the estimate is to check the nearly common estimates for
different $\alpha$ and $\phi$-divergences. Figure \ref{Co1}
presents influence functions
$\mathrm{IF}(x;T_{\alpha},P_{\theta_{0}})$, when
$\gamma\in\{-1,0,1,2,3,\}$, $\theta_{0}=0.5$ and $\alpha=0.8$. For
these choices of $\theta_{0}$ and $\alpha$, the efficiency of
$\widehat{\theta}_{n}(\alpha)$ is 0.9775 when $\gamma=1$, 0.9208
when $\gamma=2$, 0.8508 when $\gamma=3$. Here, when $\gamma$
increases, the decrease of the efficiency is worsened by a loss in
B-robustness. In this respect, the maximum likelihood estimator
appears as a good choice in terms of robustness and efficiency.

\begin{figure}[]
\begin{center}
\includegraphics[width=8cm,height=8cm]{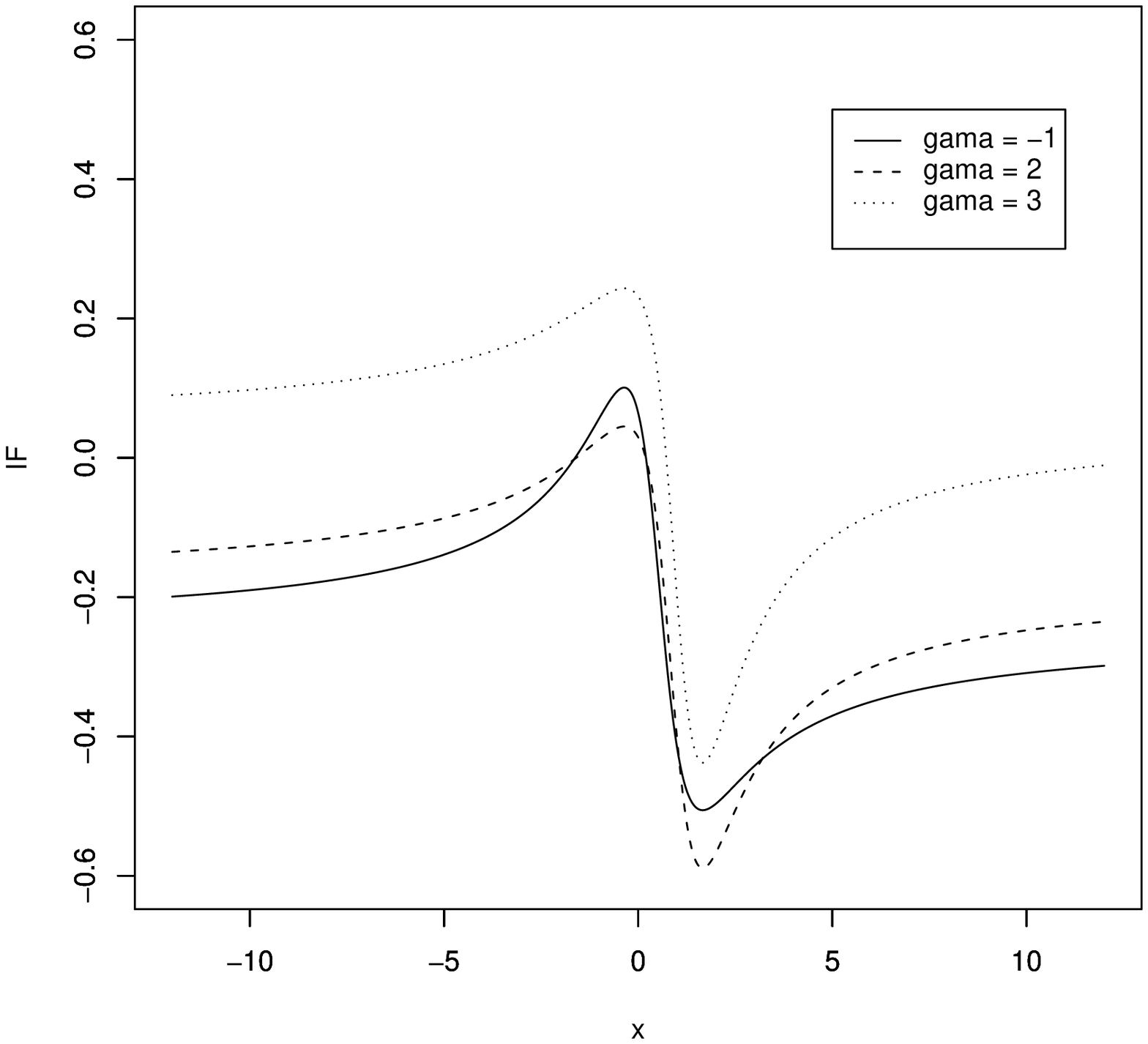}
\end{center}
\caption{Influence functions
$\mathrm{IF}(x;U_{\alpha},P_{\theta_{0}})$ for the Cauchy location
model, when the true location parameter is $\theta_{0}=0.5$ and
$\alpha=0.8$. }\label{Co2}
\end{figure}
\begin{figure}[]
\begin{center}
\includegraphics[width=8cm,height=8cm]{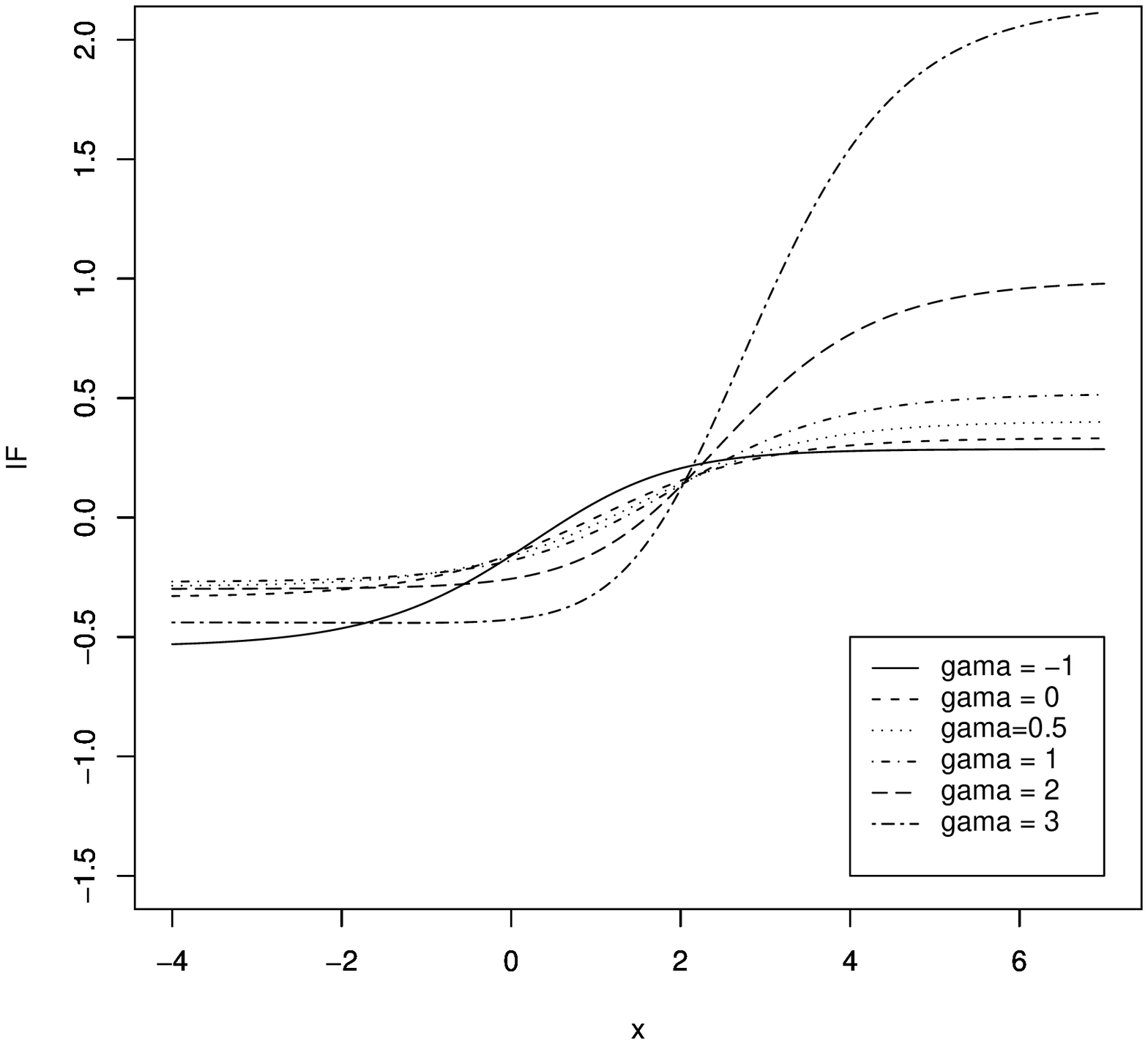}
\end{center}
\caption{Influence functions
$\mathrm{IF}(x;T_{\alpha},P_{\theta_{0}})$ for the logistic
location model, when the true location parameter is $\theta_{0}=1$
and $\alpha=1.5$. }\label{Lo1}
\end{figure}
\begin{figure}[]
\begin{center}
\includegraphics[width=8cm,height=8cm]{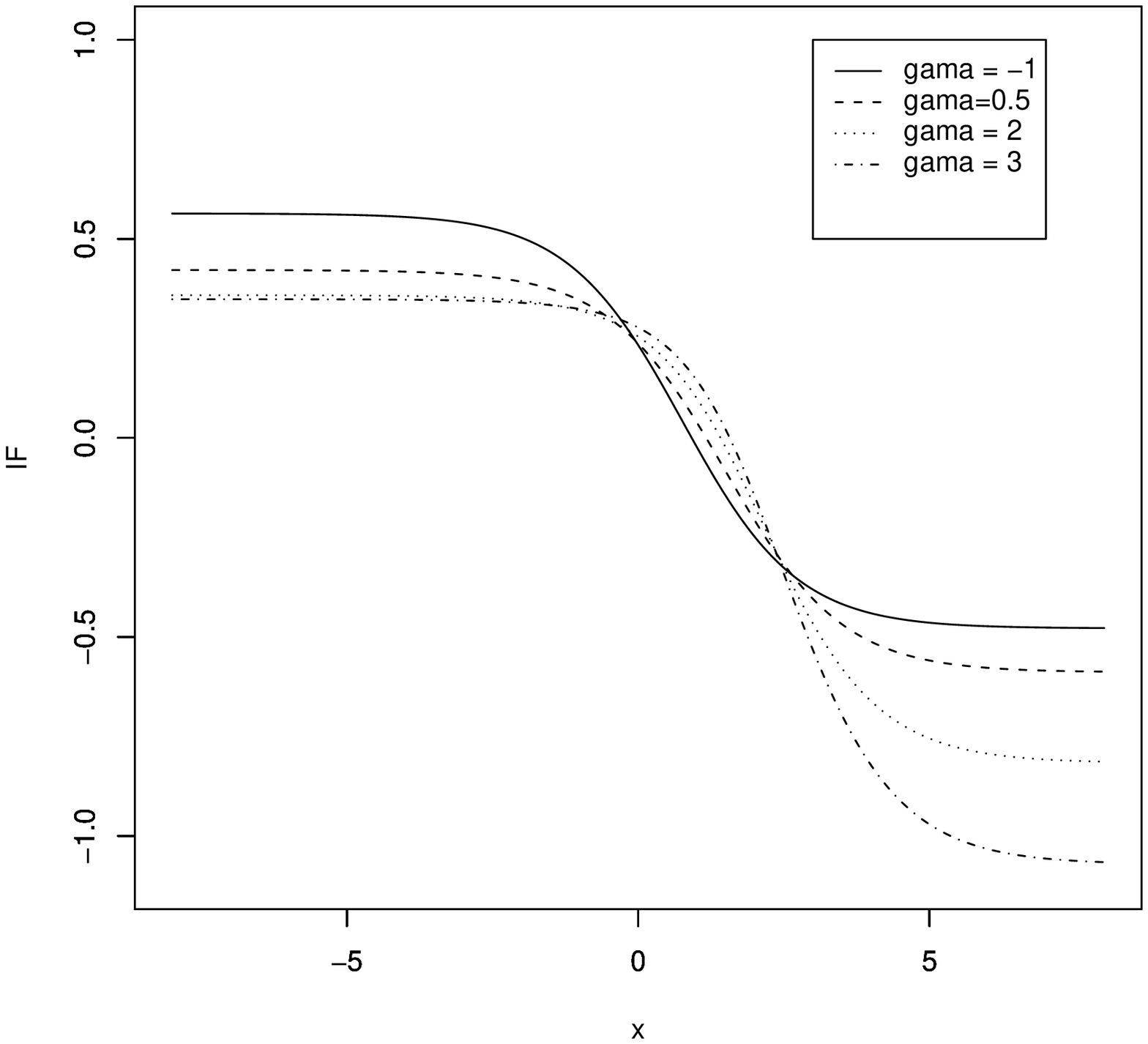}
\end{center}
\caption{Influence functions
$\mathrm{IF}(x;U_{\alpha},P_{\theta_{0}})$ for the logistic
location model, when the true location parameter is $\theta_{0}=1$
and $\alpha=1.5$. }\label{Lo2}
\end{figure}

In the case of the logistic location model, a simple calculation shows that the condition (\ref{loccond}) is fulfilled for any $\gamma$ and any $\alpha$. Also, the integral $\int m^{\prime\prime}(\theta_{0},\alpha)dP_{\theta_{0}}$ exists and is different from zero for any $\gamma$ and any $\alpha$. These conditions entail the fact that all the estimators $\widehat{\theta}_{n}(\alpha)$ are B-robust. Figure \ref{Lo1} presents influence functions $\mathrm{IF}(x;T_{\alpha},P_{\theta_{0}})$, when $\gamma\in\{-1,0,0.5,1,2,3,\}$, $\theta_{0}=1$ and $\alpha=1.5$. As in the case of the Cauchy model, when $\gamma$ increases, the decrease of the efficiency is worsened by the increase of the gross error sensitivity, such that the maximum likelihood estimator appears again as a good choice in terms of robustness and efficiency.

On the other hand, for the mean of the normal law, none of the
estimators $\widehat{\theta}_{n}(\alpha)$ is B-robust, their
influence functions being always unbounded.

In the case of the Cauchy model, as well as in the case of the logistic model, $\mathrm{IF}(x;U_{\alpha},P_{\theta_{0}})$ is bounded for any $\gamma$ and any $\alpha$. In Figure \ref{Co2}, respectively in Figure \ref{Lo2}, we present such influence functions for different choices of $\gamma$. Thus, for these two location models, all the estimators $\widehat{\phi}_{n}(\alpha,\theta_{0})$ are B-robust.

\section{Robust tests based on divergence estimators}

\subsection{Asymptotic results for contaminated observations}This subsection presents some asymptotic results
that are necessary in order to analyze the robustness of some tests based on divergence estimators. These asymptotic
results are obtained for contaminated observations, namely $X_{1},\dots,X_{n}$ are i.i.d. with
\begin{equation}\label{ccont}P_{n,\varepsilon,x}^{P}:=\left(1-\frac{\varepsilon}{\sqrt{n}}\right)P_{\theta_{n}}+\frac{\varepsilon}{\sqrt{n}}\delta_{x}\end{equation}
where $\theta_{n}=\theta_{0}+\frac{\Delta}{\sqrt{n}}$, $\Delta$
being an arbitrary vector from $\mathbb{R}^{d}$.

For $\alpha$ fixed consider the following conditions:

(C.1) The function $\theta\mapsto m(\theta,\alpha)$ is $C^{3}$ for
all $x$ and all partial derivatives of order 3 of $\theta\mapsto
m(\theta,\alpha)$ are dominated by some
$P_{\theta_{n}}$-integrable function $x\mapsto H(x)$ with the
property $\int H^{2}dP_{\theta_{n}}$ is finite, for any $n$ and
any $\Delta$.

(C.2) $\int m(\theta_{0},\alpha)dP_{\theta_{n}}$ and $\int
m^{2}(\theta_{0},\alpha)dP_{\theta_{n}}$ are finite, for any $n$
and any $\Delta$.

(C.3) $\int m^{\prime}(\theta_{0},\alpha)dP_{\theta_{n}}$ and
$\int
m^{\prime}(\theta_{0},\alpha)m^{\prime}(\theta_{0},\alpha)^{t}dP_{\theta_{n}}$
exist, for any $n$ and any $\Delta$.

(C.4) $\int m^{\prime\prime}(\theta_{0},\alpha)dP_{\theta_{n}}$
and $\int m^{\prime\prime}(\theta_{0},\alpha)^{2}dP_{\theta_{n}}$
exist, for any $n$ and any $\Delta$.

The estimators $\widehat{\theta}_{n}(\alpha)$ have good properties with respect to contamination in terms of consistency.

\begin{proposition}\label{proppas0}
If the conditions $\mathrm{(C.1), (C.3)}$ and $\mathrm{(C.4)}$ are
satisfied, then
\begin{equation*}\sqrt{n}(\widehat{\theta}_{n}(\alpha)-T_{\alpha}(P_{n,\varepsilon,x}^{P}))=O_{P}(1).\end{equation*}
\end{proposition}

Also, $\widehat{\phi}_{n}(\alpha,\theta_{0})$ enjoys normal convergence under (\ref{ccont}).

\begin{proposition}\label{proppas3}
If $\alpha\neq \theta_{0}$ and the conditions
$\mathrm{(C.1)-(C.4)}$ are satisfied, then
\begin{equation*}
\frac{\sqrt{n}(\widehat{\phi}_{n}(\alpha,\theta_{0})-U_{\alpha}(P_{n,\varepsilon,x}^{P}))}{[\int
\mathrm{IF}^{2}(y;U_{\alpha},P_{n,\varepsilon,x}^{P})dP_{n,\varepsilon,x}^{P}(y)]^{1/2}}
\end{equation*}
converges in distribution to a normal standard variable.
\end{proposition}
\subsection{Robust tests based on divergence estimators}
In this subsection we propose tests based on dual
$\phi$-divergence estimators and study their robustness
properties. We mention that the use of the dual form of a
divergence to derive robust tests was discussed in a different
context by Broniatowski and Leorato \cite{BrLe06} in the case of
the Neyman $\chi^{2}$ divergence.

For testing the hypothesis $\theta=\theta_{0}$ against the
alternative $\theta\neq\theta_{0}$, consider the test of level
$\alpha_{0}$ defined by the test statistic
$\widehat{\phi}_{n}:=\widehat{\phi}_{n}(\alpha,\theta_{0})$ with
$\alpha\neq \theta_{0}$ and by the critical region
\begin{equation*}
C:=\left\{\left|\frac{\sqrt{n}(\widehat{\phi}_{n}-\phi(\alpha,\theta_{0}))}{[\int\mathrm{IF}^{2}(y;U_{\alpha},P_{\theta_{0}})dP_{\theta_{0}}(y)]^{1/2}}\right|\geq
q_{1-\frac{\alpha_{0}}{2}} \right\}
\end{equation*}
where $q_{1-\frac{\alpha_{0}}{2}}$ is the
$(1-\frac{\alpha_{0}}{2})$-quantile of the standard normal
distribution.

Due to the asymptotic normality of $\widehat{\phi}_{n}$, for $n$
large, the level writes as
\begin{eqnarray}
\alpha_{0}&\simeq&P_{\theta_{0}}\left(\left|\frac{\sqrt{n}(\widehat{\phi}_{n}-\phi(\alpha,\theta_{0}))}{[\int\mathrm{IF}^{2}(y;U_{\alpha},P_{\theta_{0}})dP_{\theta_{0}}(y)]^{1/2}}\right|\geq
q_{1-\frac{\alpha_{0}}{2}}\right)\\ &=&
P_{\theta_{0}}(|\widehat{\phi}_{n}-\phi(\alpha,\theta_{0})|\geq
(\sqrt{n})^{-1}[\int\mathrm{IF}^{2}(y;U_{\alpha},P_{\theta_{0}})dP_{\theta_{0}}(y)]^{1/2}q_{1-\frac{\alpha_{0}}{2}})\\
&=&2P_{\theta_{0}}(\widehat{\phi}_{n}\geq
k_{n}(\alpha_{0}))\label{llevel}
\end{eqnarray}
where
$k_{n}(\alpha_{0})=(\sqrt{n})^{-1}[\int\mathrm{IF}^{2}(y;U_{\alpha},P_{\theta_{0}})dP_{\theta_{0}}(y)]^{1/2}q_{1-\frac{\alpha_{0}}{2}}+\phi(\alpha,\theta_{0})$.

We work with the form (\ref{llevel}) of the level and consequently
of the probability to reject the null hypothesis, this being
easier to handle in the proofs of the results that follows.

Consider the sequence of contiguous alternatives
$\theta_{n}=\theta_{0}+\Delta n^{-1/2}$, where $\Delta$ is any
vector from $\mathbb{R}^{d}$. When $\theta_{n}$ tends to
$\theta_{0}$, the contamination must converge to 0 at the same
rate, to avoid the overlapping between the neighborhood of the
hypothesis and that of the alternative (see Hampel et al.
\cite{HaRoRoSt}, p.198 and Heritier and Ronchetti \cite{HeRo}).
Therefore we consider the contaminated distributions
\begin{equation}\label{ld}
P_{n,\varepsilon,x}^{L}=\left(1-\frac{\varepsilon}{\sqrt{n}}\right)P_{\theta_{0}}+\frac{\varepsilon}{\sqrt{n}}\delta_{x}
\end{equation}
for the level and
\begin{equation}\label{pd}
P_{n,\varepsilon,x}^{P}=\left(1-\frac{\varepsilon}{\sqrt{n}}\right)P_{\theta_{n}}+\frac{\varepsilon}{\sqrt{n}}\delta_{x}
\end{equation}
for the power.

The asymptotic level (the asymptotic power) under (\ref{ld}) (under (\ref{pd})) will be evaluated now.

Let $\beta_{0}=\lim_{n\to
\infty}2P_{\theta_{n}}(\widehat{\phi}_{n}\geq k_{n}(\alpha_{0}))$
be the asymptotic power of the test under the family of
alternatives $P_{\theta_{n}}$. The test is robust with respect to
the power if the limit of the powers under the contaminated
alternatives stays in a bounded neighborhood of $\beta_{0}$, so
that the role of the contamination is somehow controlled. Also,
the test is robust with respect to the level if the limit of the
level under the contaminated null distributions stays in a bounded
neighborhood of $\alpha_{0}$.

Let
$P_{n,\varepsilon,x}=2P_{n,\varepsilon,x}^{P}(\widehat{\phi}_{n}\geq
k_{n}(\alpha_{0}))$. In the same vein as in Dell'Aquilla and
Ronchetti \cite{DeRo} it holds:
\begin{proposition}\label{p11}
If the conditions $\mathrm{(C.1)-(C.4)}$ are fulfilled, then the
asymptotic power of the test under $P_{n,\varepsilon,x}^{P}$ is
given by
\begin{eqnarray}
\label{t1}\lim_{n\to
\infty}P_{n,\varepsilon,x}&=&2-2\Phi\left(\Phi^{-1}\left(1-\frac{\alpha_{0}}{2}\right)-\Delta\frac{c}{[\int
\mathrm{IF}^{2}(y;U_{\alpha},P_{\theta_{0}})dP_{\theta_{0}}(y)
]^{1/2}}-\right.\nonumber\\&&-\left.\varepsilon\frac{\mathrm{IF}(x;U_{\alpha},P_{\theta_{0}})}{[\int
\mathrm{IF}^{2}(y;U_{\alpha},P_{\theta_{0}})dP_{\theta_{0}}(y)\label{asp}
]^{1/2}}\right)
\end{eqnarray}
where $c=\int
m(\theta_{0},\alpha,y)\frac{\dot{p}_{\theta_{0}}(y)}{p_{\theta_{0}}(y)}dP_{\theta_{0}}(y)$
and $\Phi$ is the cumulative distribution function of the standard
normal.
\end{proposition}

A Taylor expansion with respect to $\varepsilon$ yields
\begin{eqnarray*}
&&\!\!\!\!\!\!\!\!\!\!\!\!\lim_{n\to \infty}P_{n,\varepsilon,x} =
2-2\Phi\left(\Phi^{-1}\left(1-\frac{\alpha_{0}}{2}\right)-\Delta\frac{c}{[\int
\mathrm{IF}^{2}(y;U_{\alpha},P_{\theta_{0}})dP_{\theta_{0}}(y)
]^{1/2}}\right)+ \\  &&\!\!\!\!\!\!\!\!\!\!\!\!+2\varepsilon
f\left(\!\Phi^{-1}\left(1\!-\!\frac{\alpha_{0}}{2}\right)\!-\!\Delta\frac{c}{[\int
\mathrm{IF}^{2}(y;U_{\alpha},P_{\theta_{0}})dP_{\theta_{0}}(y)
]^{1/2}}\!\right)\!\!\frac{\mathrm{IF}(x;U_{\alpha},P_{\theta_{0}}\!)}{[\int
\mathrm{IF}^{2}(y;U_{\alpha},P_{\theta_{0}})dP_{\theta_{0}}(y)
]^{1/2}}\!+\!o(\varepsilon)\\
&&\!\!\!\!\!\!\!\!\!\!\!\!=\beta_{0}\!+\!2\varepsilon
f\left(\Phi^{-1}\left(1\!-\!\frac{\alpha_{0}}{2}\!\right)\!-\!\Delta\frac{c}{[\int
\mathrm{IF}^{2}(y;U_{\alpha},P_{\theta_{0}})dP_{\theta_{0}}(y)
]^{1/2}}\right)\frac{\mathrm{IF}(x;U_{\alpha},P_{\theta_{0}})}{[\int
\mathrm{IF}^{2}(y;U_{\alpha},P_{\theta_{0}})dP_{\theta_{0}}(y)
]^{1/2}}+\\&&+o(\varepsilon)
\end{eqnarray*}
where $\beta_{0}$ is the asymptotic power for the non contaminated model and $f$ is
the density of the standard normal distribution.

In order to limit the
bias in the power of the test it is sufficient to bound the influence function
$\mathrm{IF}(x;U_{\alpha}, P_{\theta_{0}})$. Bounding the influence function is
therefore enough to maintain the power in a pre-specified band
around $\beta_{0}$.

Let
$L_{n,\varepsilon,x}=2P_{n,\varepsilon,x}^{L}(\widehat{\phi}_{n}\geq
k_{n}(\alpha_{0}))$. Putting $\Delta=0$ in (\ref{t1}) yields:
\begin{proposition}\label{p12}
If the conditions $\mathrm{(C.1)-(C.4)}$ are fulfilled, then the
asymptotic level of the test under $P_{n,\varepsilon,x}^{L}$ is
given by
\begin{eqnarray*}
\lim_{n\to
\infty}L_{n,\varepsilon,x}&=&2-2\Phi\left(\Phi^{-1}\left(1-\frac{\alpha_{0}}{2}\right)-\varepsilon\frac{\mathrm{IF}(x;U_{\alpha},P_{\theta_{0}})}{[\int
\mathrm{IF}^{2}(y;U_{\alpha},P_{\theta_{0}})dP_{\theta_{0}}(y)
]^{1/2}}\right)\\ &=&\alpha_{0}+\varepsilon
f\left(\Phi^{-1}\left(1-\frac{\alpha_{0}}{2}\right)\right)\frac{\mathrm{IF}(x;U_{\alpha},P_{\theta_{0}})}{[\int
\mathrm{IF}^{2}(y;U_{\alpha},P_{\theta_{0}})dP_{\theta_{0}}(y)]^{1/2}}+o(\varepsilon).
\end{eqnarray*}
\end{proposition}
Hence, when $\mathrm{IF}(x;U_{\alpha},P_{\theta_{0}})$ is bounded, $L_{n,\varepsilon,x}$ remains
between pre-specified bounds of $\alpha_{0}$.

As the Proposition \ref{p11} and Proposition \ref{p12} show, both
the asymptotic power of the test under $P_{n,\varepsilon,x}^{P}$
and the asymptotic level of the test under
$P_{n,\varepsilon,x}^{L}$ are controlled by the influence function
of the test statistic. Hence, the robustness of the test statistic
$\widehat{\phi}_{n}$, as discussed in the previous section,
assures the stability of the test under small arbitrary departures
from the null hypothesis, as well as a good power under small
arbitrary departures from specified alternatives. Figures
\ref{Co2} and \ref{Lo2} provide some specific values of $\gamma$
and $\alpha$ inducing robust tests for $\theta_{0}$ corresponding
to those models.

\section{Simulation results}Simulation were run in order to
examine empirically the performances of the robust dual
$\phi$-divergence estimators and tests. The considered parametric
model was the scale normal model with known mean. We worked with
data generated from the model, as well as with contaminated data.

To make some comparisons, beside dual $\phi$-divergence
estimators, we considered minimum density power divergence
estimators of Basu et al. \cite{BaHa} (MDPDE's) and the maximum
likelihood estimator (MLE). Recall that a MDPDE of a parameter
$\theta$ is obtained as solution of the equation
\begin{equation}\label{Baseq}
\int
\dot{p}_{\theta}(z)p_{\theta}^{\beta}(z)dz-\frac{1}{n}\sum_{i=1}^{n}\dot{p}_{\theta}(X_{i})p_{\theta}^{\beta-1}(X_{i})=0
\end{equation}
with respect to $\theta$, where $\beta>0$ and $X_{1},\dots,X_{n}$
is a sample from $P_{\theta}$. In the case of the scale normal
model $\mathcal{N}(m,\sigma)$, equation (\ref{Baseq}) writes
as\begin{eqnarray*} \hspace{5mm}&&\int
\frac{1}{\sigma^{\beta+2}(\sqrt{2\pi})^{\beta+1}}\left(e^{-\frac{1}{2}\left(\frac{z-m}{\sigma}\right)^{2}}\right)^{\beta+1}\left[-1+\left(\frac{z-m}{\sigma}\right)^{2}\right]dz-\\
\hspace{5mm}&&-\frac{1}{n}\sum_{i=1}^{n}\frac{1}{\sigma^{\beta+1}(\sqrt{2\pi})^{\beta}}\left(e^{-\frac{1}{2}\left(\frac{X_{i}-m}{\sigma}\right)^{2}}\right)^{\beta}\left[-1+\left(\frac{X_{i}-m}{\sigma}\right)^{2}\right]=0
\end{eqnarray*}
and the MDPDE of the parameter $\sigma$ is robust for any
$\beta>0$.

In a first Monte Carlo experiment the data were generated from the
scale normal model $\mathcal{N}(0,1)$ with mean $m=0$ known,
$\sigma=1$ being the parameter of interest. We considered
different choices for the tuning parameter $\alpha$ and for the
Cressie-Read divergence to compute D$\phi$E's, and different
choices for the tuning parameter $\beta$ in order to compute
MDPDE's. For each set of configurations considered, 5000 samples
of size $n=100$ were generated from the model, and for each sample
D$\phi$E's, MDPDE's and MLE were obtained.

In Table 1 we present the results of the simulations, showing
simulation based estimates of the bias and MSE given by
\begin{equation*}
\widehat{\mathrm{Bias}}=\frac{1}{n_{s}}\sum_{i=1}^{n_{s}}(\widehat{\sigma}_{i}-\sigma),\;\;\widehat{\mathrm{MSE}}=\frac{1}{n_{s}}\sum_{i=1}^{n_{s}}(\widehat{\sigma}_{i}-\sigma)^{2},
\end{equation*}
where $n_{s}$ denotes the number of samples (5000 in our case) and
$\widehat{\sigma}_{i}$ denotes an estimate of $\sigma$ for the
$i$th sample. Examination of the table shows that D$\phi$E's give
as good results as MDPDE's or MLE.

In a second Monte Carlo experiment, we first generated samples
with 100 observations, namely 98 coming from $\mathcal{N}(0,1)$
and 2 outliers $x=10$ and then we generated samples with 100
observations, namely 96 from $\mathcal{N}(0,1)$ and 4 outliers
$x=10$. The tuning parameters were the same as in the non
contaminated case and also $n_{s}=5000$. The simulation results
are given in Table 2. As can be seen, the results for D$\phi$E's
and MDPDE's are comparable, they being better than the results for
MLE in both cases.

A close look at the results of the simulations show the D$\phi$E
performs well under the model, when no outliers are generated;
indeed the best results are obtained when $\gamma=-0.1$, whatever
$\overline{\sigma}=1.5$ or $\overline {\sigma}=1.9$. The
performance of the estimator under the model is comparable to that
of some MDPDE's in terms of empirical MSE ($\mathrm{\widehat
{MSE}}$): indeed the $\mathrm{\widehat{MSE}}$ for D$\phi$E with
$\gamma=-0.1$ parallels MDPDE's for small $\beta.$ It is also
slightly shorter than the one obtained through the MLE. Under
contamination, the D$\phi$E with $\gamma=-0.5$ yields clearly the
most robust estimate and the empirical MSE is very small,
indicating a strong stability of the estimate. It compares
favorably with MDPE for all $\beta,$ whatever
$\overline{\sigma}=1.5$ or $\overline{\sigma}=1.9$. The simulation
with $4$ outliers at $x=10$ provide a clear evidence of the
properties of the D$\phi$E with $\gamma=-0.5$. Also small values
of $\beta$ give similar results as large negative values of
$\gamma$, whatever $\overline{\sigma}$, under contamination.
Although $\gamma=-0.1$ is a good alternative to MLE under the
model, $\gamma=-0.5$ behaves quite well in terms of bias while
keeping short empirical MSE under the model or under
contamination. These results are in full accordance with Figure 1;
indeed the influence function is constant close to $0$ for large
values of $x$.

Thus, the D$\phi$E is shown to be an attractive alternative to
both the MLE and MDPDE in these settings.

In order to test the hypothesis $\sigma=1$ with respect to the
alternative $\sigma\neq 1$, we considered the test statistic
\begin{equation*}
\frac{\sqrt{n}(\widehat{\phi}_{n}-\phi(\alpha,\theta_{0}))}{[\int\mathrm{IF}^{2}(y;U_{\alpha},P_{\theta_{0}})dP_{\theta_{0}}(y)]^{1/2}}
\end{equation*}
(here $\theta_{0}=\sigma=1$). Under the null hypothesis, this test
statistic is asymptotically $\mathcal{N}(0,1)$. We worked with
data generated from the model $\mathcal{N}(0,1)$, as well as with
contaminated data. In each case, we simulated 5000 samples and we
computed the actual levels
\begin{equation*}
P\left(\left|\frac{\sqrt{n}(\widehat{\phi}_{n}-\phi(\alpha,\theta_{0}))}{[\int\mathrm{IF}^{2}(y;U_{\alpha},P_{\theta_{0}})dP_{\theta_{0}}(y)]^{1/2}}\right|\geq
q_{1-\frac{\alpha_{0}}{2}}\right)
\end{equation*}
corresponding to the nominal levels
$\alpha_{0}=0.01,0.02,\dots,0.1$. We reported the corresponding
relative errors
\begin{equation*}
\left(P\left(\left|\frac{\sqrt{n}(\widehat{\phi}_{n}-\phi(\alpha,\theta_{0}))}{[\int\mathrm{IF}^{2}(y;U_{\alpha},P_{\theta_{0}})dP_{\theta_{0}}(y)]^{1/2}}\right|\geq
q_{1-\frac{\alpha_{0}}{2}}\right)-\alpha_{0}\right)/\alpha_{0}.
\end{equation*}
In Figure \ref{Gr7} we present relative errors for the robust
tests applied to the scale normal model $\mathcal{N}(0,1)$, when
the data are generated from the model. The sample size is $n=100$,
the tuning parameter is $\overline{\sigma}=1.9$ and the
Cressie-Read divergences correspond to
$\gamma\in\{-1.5,-1,-0.5,-0.1\}$. The approximation of the level
is good for all the considered divergences.

In Figure \ref{Gr8} are represented relative errors of the robust
tests applied to the scale normal model $\mathcal{N}(0,1)$, for
samples with $n=100$ data, namely 98 data generated from
$\mathcal{N}(0,1)$ and 2 outliers $x=10$. We considered
$\overline{\sigma}=1.9$ and $\gamma\in\{-2,-1.5\}$. Again, the
approximation of the level of the test is good for all the
considered divergences.

In Figure \ref{Gr9} we present relative errors of the robust tests
applied to the scale normal model $\mathcal{N}(0,1)$, for samples
with $n=100$ data, namely 96 data generated from
$\mathcal{N}(0,1)$ and 4 outliers $x=10$. We considered
$\overline{\sigma}=1.9$ and $\gamma\in\{-2,-1.5\}$.

Observe that the tests give good results for values of $\gamma$
close to zero when the data are not contaminated, respectively for
large negative values of $\gamma$ when the data are contaminated.

Thus, the numerical results show that dual $\phi$-divergence
estimates and corresponding tests are stable in the presence of
some outliers in the sample.

\begin{table}[tbp] Table 1.\vspace{3mm}\\ \text{Simulation results for D$\phi$E, MDPDE and MLE of}\\
\text{the parameter $\sigma=1$ when the data are
generated}\\\vspace{2mm}\text{from the model
$\mathcal{N}(0,1)$.}\\
\begin{tabular}{llllll}
\hline \vspace{1mm} &  & $\widehat{\sigma}$ &
$\widehat{\mathrm{Bias}}$ & $\widehat{\mathrm{MSE}}$
\\ \hline D$\phi$E &  &  &  &
\\
$\overline{\sigma}$=1.5\;\;$\gamma=-2$ &  & 0.99770 & -0.00229
&0.00917
      \\

$\overline{\sigma}$=1.5\;\;$\gamma=-1.5$ &  & 0.99735 & -0.00264
&0.00822
\\

$\overline{\sigma}$=1.5\;\;$\gamma=-1$ &  & 0.99760 & -0.00239
&0.00698
\\

$\overline{\sigma}$=1.5\;\;$\gamma=-0.5$ &  & 0.99833 & -0.00166
&0.00563
\\

$\overline{\sigma}$=1.5\;\;$\gamma=-0.1$ &  & 0.99799 & -0.00200
&0.00492

\\ $\overline{\sigma}$=1.9\;\;$\gamma=-2$ &  & 0.99892 & -0.00107 &
0.01029      \\

$\overline{\sigma}$=1.9\;\;$\gamma=-1.5$ &  & 0.99841 & -0.00158 &
0.00924      \\

$\overline{\sigma}$=1.9\;\;$\gamma=-1$ &  & 0.99824 & -0.00175 &
0.00773     \\

$\overline{\sigma}$=1.9\;\;$\gamma=-0.5$ &  & 0.99839 & -0.00160 &
0.00588     \\

$\overline{\sigma}$=1.9\;\;$\gamma=-0.1$ &  & 0.99768 & -0.00231 &
0.00473     \\

MDPDE &  &  &  &      \\ $\beta=0.1$ & &0.99894 & -0.00105 &
0.00514
\\ $\beta=0.5$ &  & 0.99986 & -0.00013 &0.00686
   \\ $\beta=1$ &  & 1.00005 & \;0.00005 & 0.00927

\\
$\beta=1.5$ &  & 1.00074 & \;0.00074 & 0.01077
\\ $\beta=2$ &  & 1.00150 & \;0.00150& 0.01165

\\
$\beta=2.5$ &  & 1.00294 & \;0.00294 & 0.01266
\\ MLE &  & 0.99743 & -0.00256 & 0.00501    \\ \hline
\end{tabular}
\vspace{2cm}
\end{table}
\newpage

\begin{center}
\begin{table}[tbp] Table 2.\vspace{3mm}\\ \text{Simulation
results for D$\phi$E, MDPDE and MLE of the parameter $\sigma=1$
when 98 data}\\ \text{are generated from the model
$\mathcal{N}(0,1)$ and 2 outliers $x=10$ are added,
respectively}\\\vspace{2mm}\text{when 96 data are generated from
the model $\mathcal{N}(0,1)$ and 4 outliers $x=10$ are added.}\\
\begin{tabular}{lllllllll} \hline
&&&\!2 outliers\!&&&&\!4 outliers&\\ \cline{3-5}\cline{7-9}
\vspace{1mm} &  & $\widehat{\sigma}$ & $\widehat{\mathrm{Bias}}$ &
$\widehat{\mathrm{MSE}}$&& $\widehat{\sigma}$ &
$\widehat{\mathrm{Bias}}$ & $\widehat{\mathrm{MSE}}$    \\ \hline
D$\phi$E & & & &&&&&
\\

 $\overline{\sigma}$=1.5\;\;$\gamma=-2$ &  & 1.01186 & 0.01186
 &0.00914& &1.02540  & \;0.02540 & 0.00946
       \\

$\overline{\sigma}$=1.5\;\;$\gamma=-1.5$ &  & 1.00850 & 0.00850
&0.00816 &&1.01911 & \;0.01911 &0.00833
\\

$\overline{\sigma}$=1.5\;\;$\gamma=-1$ &  & 1.00499 & 0.00499
&0.00697 && 1.01210 & \;0.01210 &0.00707
\\

$\overline{\sigma}$=1.5\;\;$\gamma=-0.5$ &  & 1.00171 & 0.00171
&0.00572 && 1.00526 & \;0.00526 &0.00583
\\

$\overline{\sigma}$=1.5\;\;$\gamma=-0.1$ &  & 1.09661 & 0.09661
&0.01641 && 0.99766 & -0.00233 &0.00088  \\
$\overline{\sigma}$=1.9\;\;$\gamma=-2$ &  & 1.01589 & 0.01589 &
0.01059 &&1.03547  & \;0.03547 &0.01182        \\

$\overline{\sigma}$=1.9\;\;$\gamma=-1.5$ &  & 1.01236 & 0.01236 &
0.00942  && 1.02840 & \;0.02840 &0.01027     \\

$\overline{\sigma}$=1.9\;\;$\gamma=-1$ &  & 1.00785 & 0.00785 &
0.00785  && 1.01912 & \;0.01912 & 0.00838    \\

$\overline{\sigma}$=1.9\;\;$\gamma=-0.5$ &  & 1.00274 & 0.00274 &
0.00598   && 1.00842 & \;0.00842 &0.00637    \\

$\overline{\sigma}$=1.9\;\;$\gamma=-0.1$ &  & 1.06708 & 0.06708 &
0.02241 & & 1.10531 & \;0.10531 &0.02083     \\

MDPDE &  &  &  &   &&&&    \\ $\beta=0.1$ &  & 1.01117 & 0.01117 &
0.00646 && 1.02676 & \;0.02676& 0.00891
\\ $\beta=0.5$ &  & 1.00700 & 0.00700 & 0.00712 && 1.01417 & \;0.01417 & 0.00743
   \\ $\beta=1$ &  & 1.01406 &0.01406  & 0.00975 &&1.02892 & \;0.02892 & 0.01062

\\
$\beta=1.5$ &  & 1.01916 & 0.01916 & 0.01148 && 1.03876 &
\;0.03876 & 0.01297
\\ $\beta=2$ &  & 1.02233 &0.02233 & 0.01254 && 1.04448 & \;0.04447 & 0.01450

\\
$\beta=2.5$ &  & 1.02450 & 0.02450 & 0.01342 && 1.04771 &\;0.04771
& 0.01556
\\ MLE &  & 1.72587 & 0.72587 & 0.52852 && 2.22720 & \;1.22720 & 1.50701  \\ \hline\end{tabular}
\vspace{2cm}
\end{table}
\end{center}

\newpage
\begin{center}
{\scriptsize
\begin{figure}[tbp]
{\scriptsize
\begin{tabular}{cc}
\includegraphics[width=6.5cm,height=6cm]{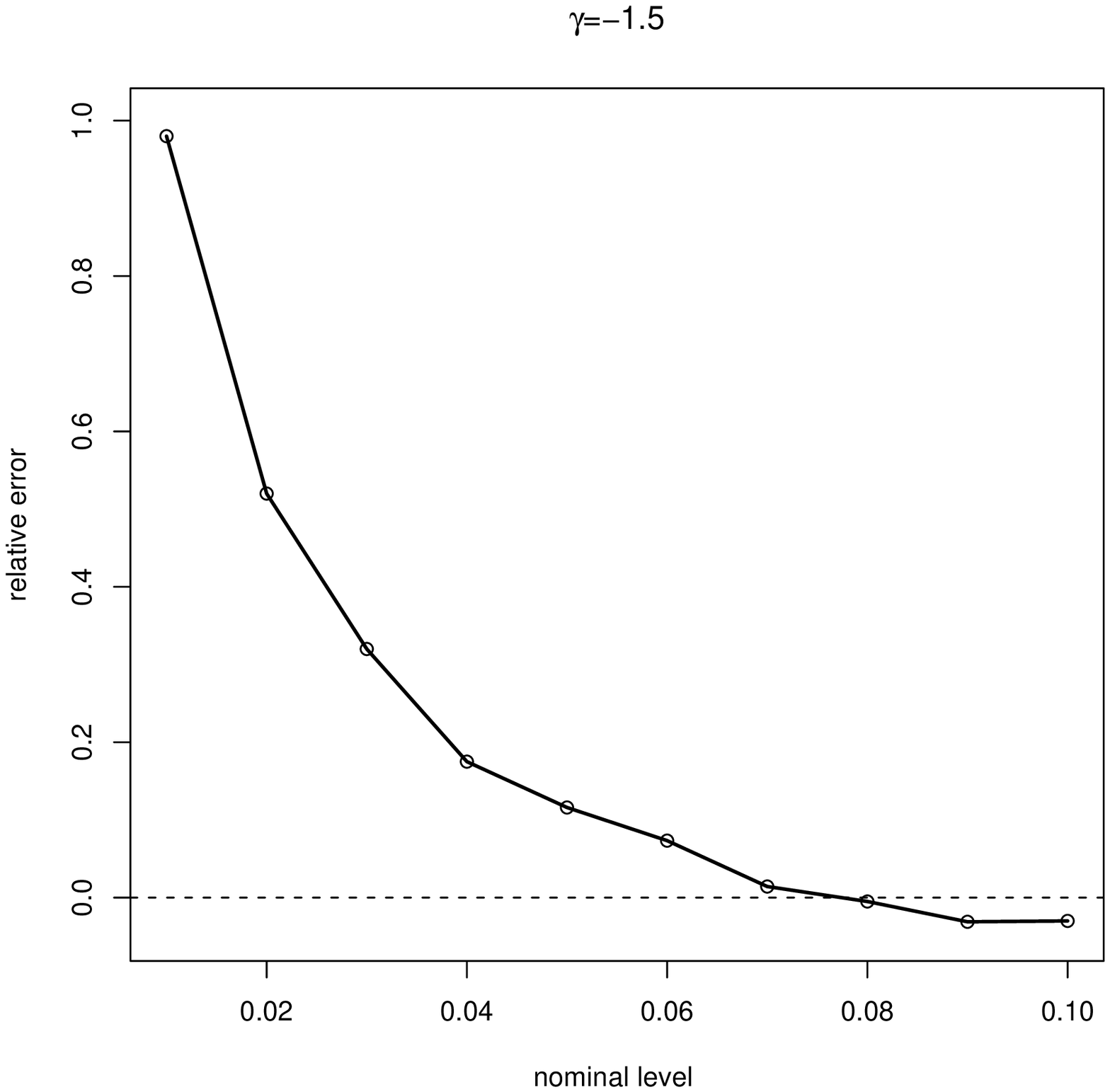} &
\includegraphics[width=6.5cm,height=6cm]{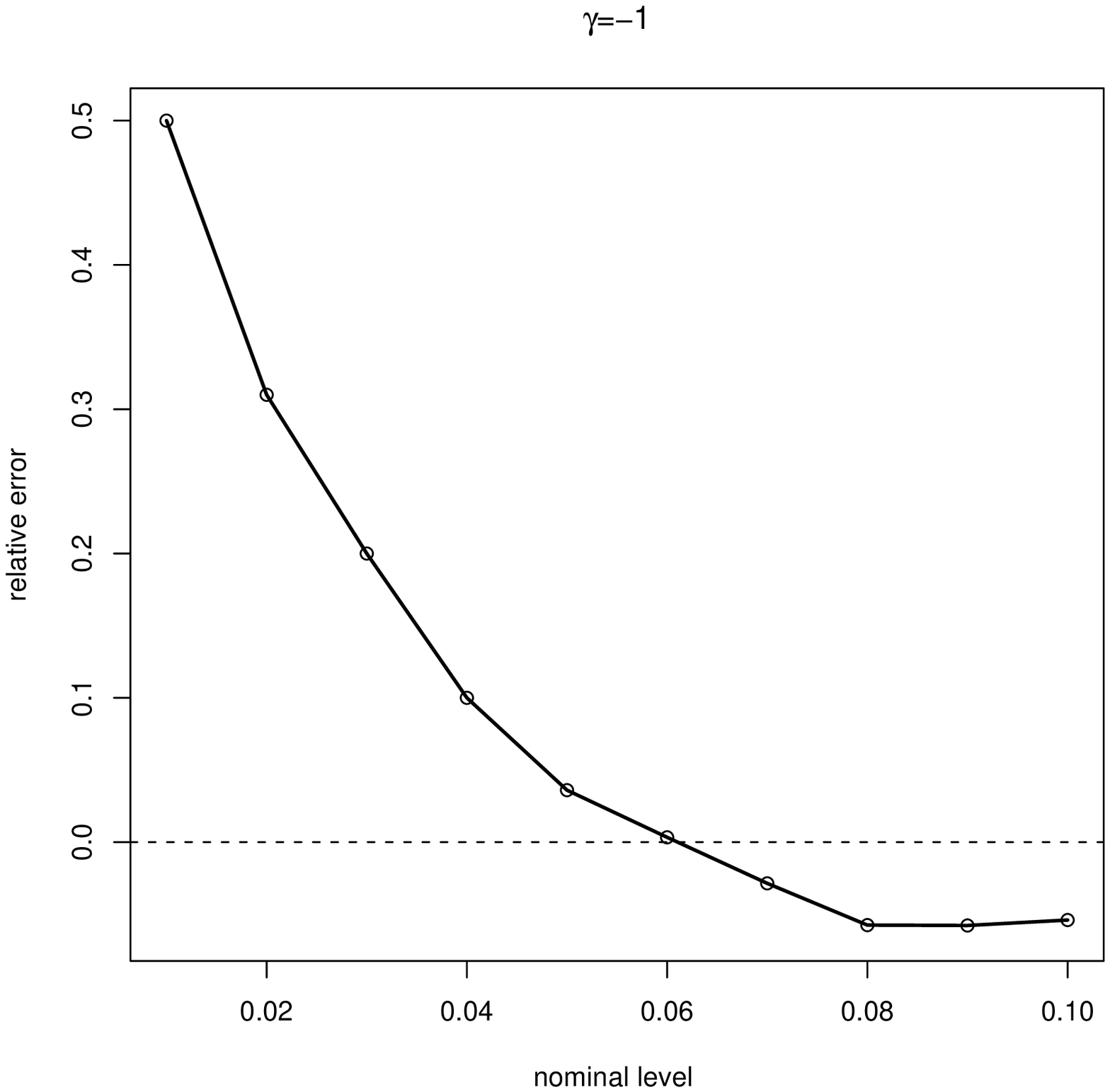} \\
\includegraphics[width=6.5cm,height=6cm]{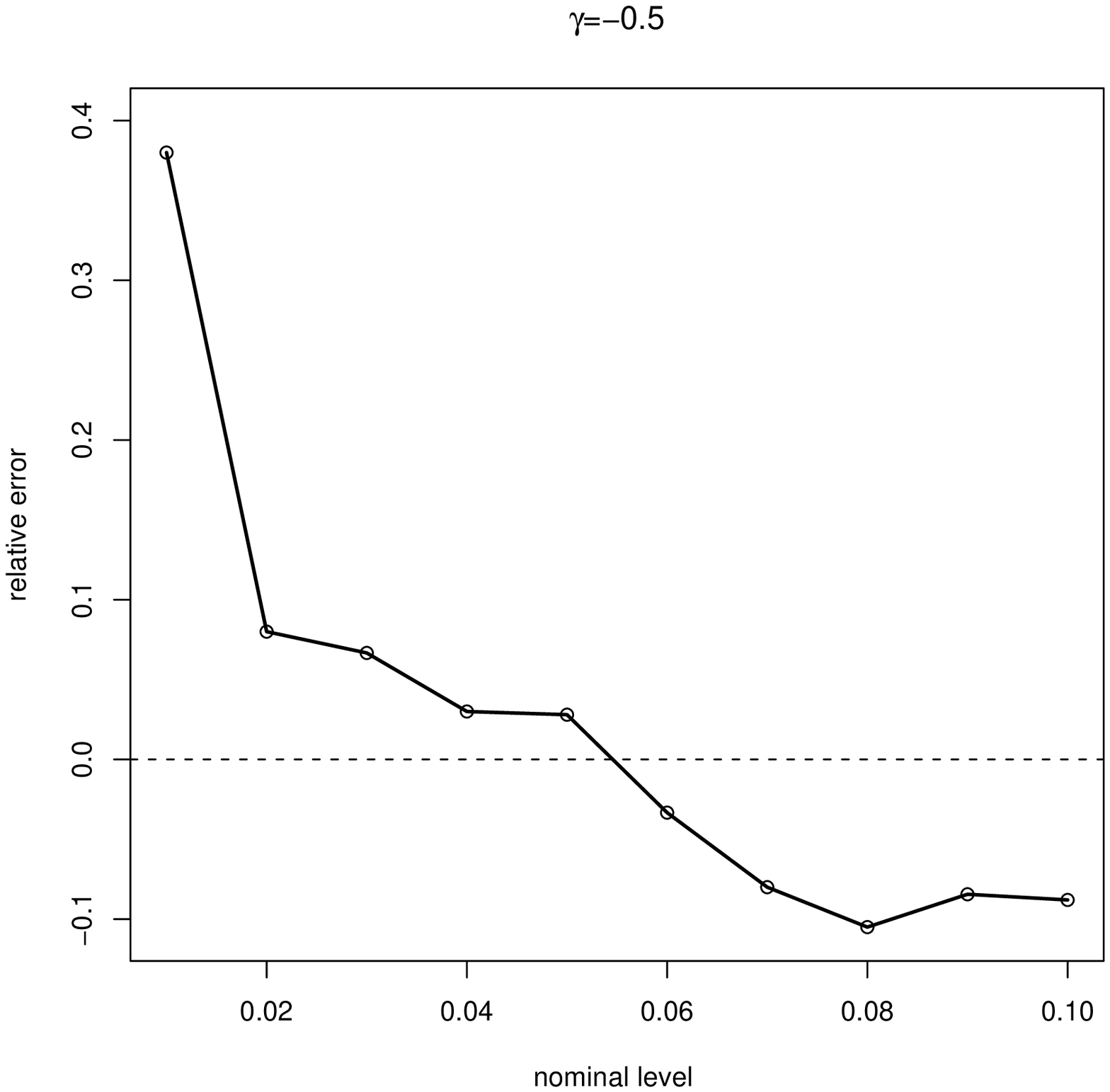} &
\includegraphics[width=6.5cm,height=6cm]{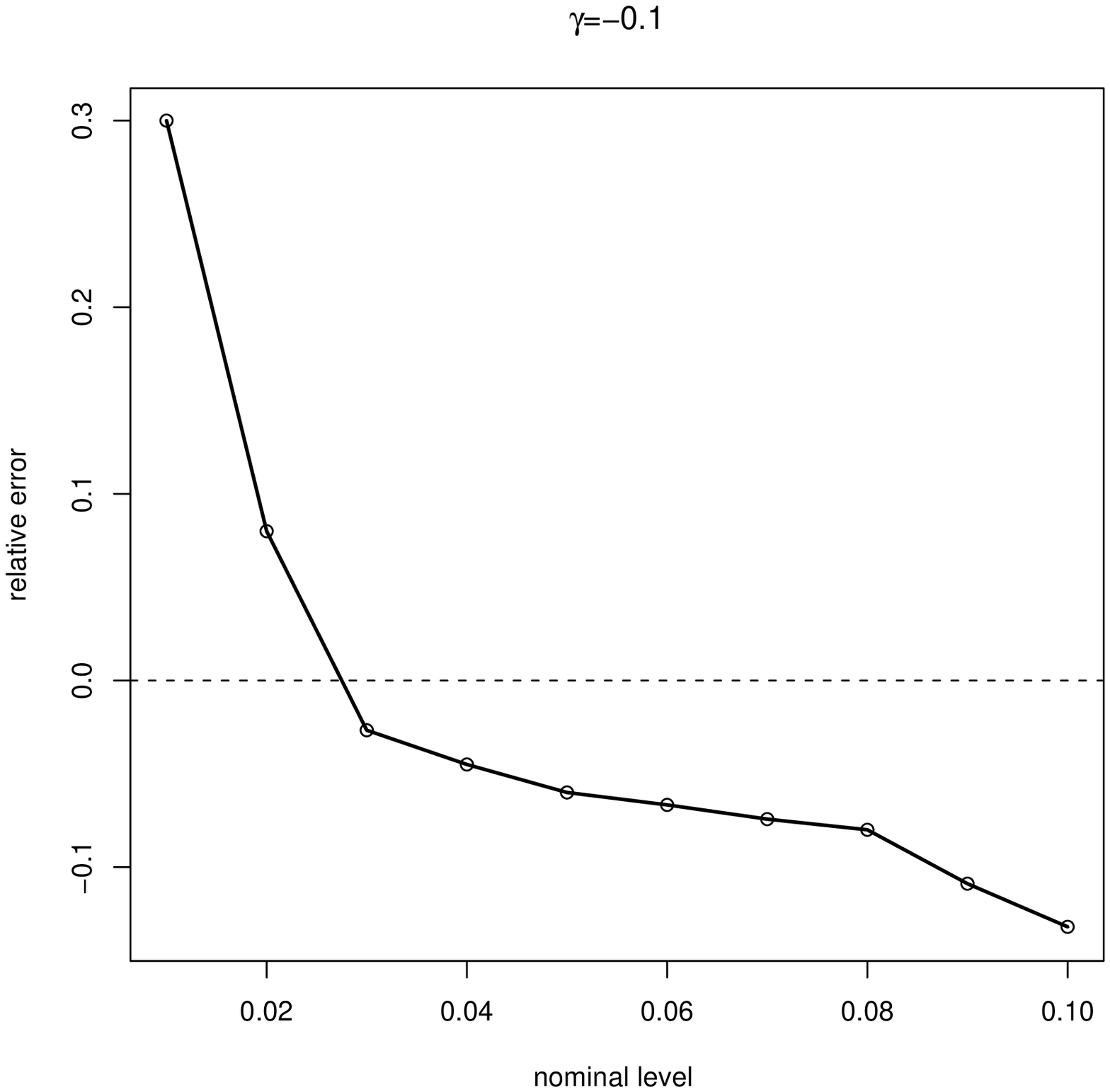}
\end{tabular}
\vspace{-2mm}} \caption{Relative errors of the robust tests
applied to the scale normal model $\mathcal{N}(0,1)$, when
$\overline{\sigma}=1.9$ and 100 data are generated from
model.}\label{Gr7}
\end{figure}
}
\end{center}

\begin{center}
{\scriptsize
\begin{figure}[tbp]
{\scriptsize
\begin{tabular}{cc}
\includegraphics[width=6.5cm,height=6cm]{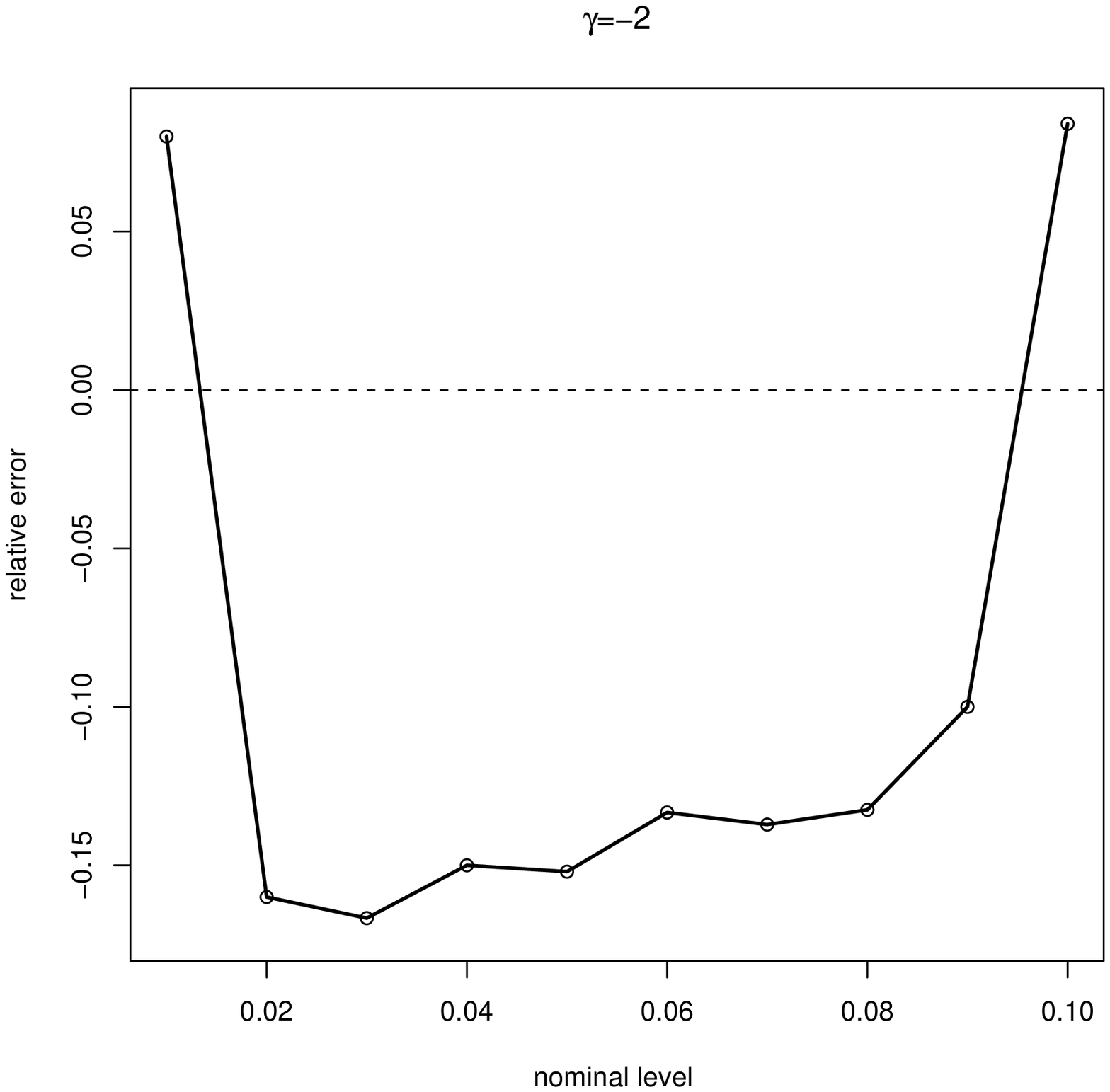} &
\includegraphics[width=6.5cm,height=6cm]{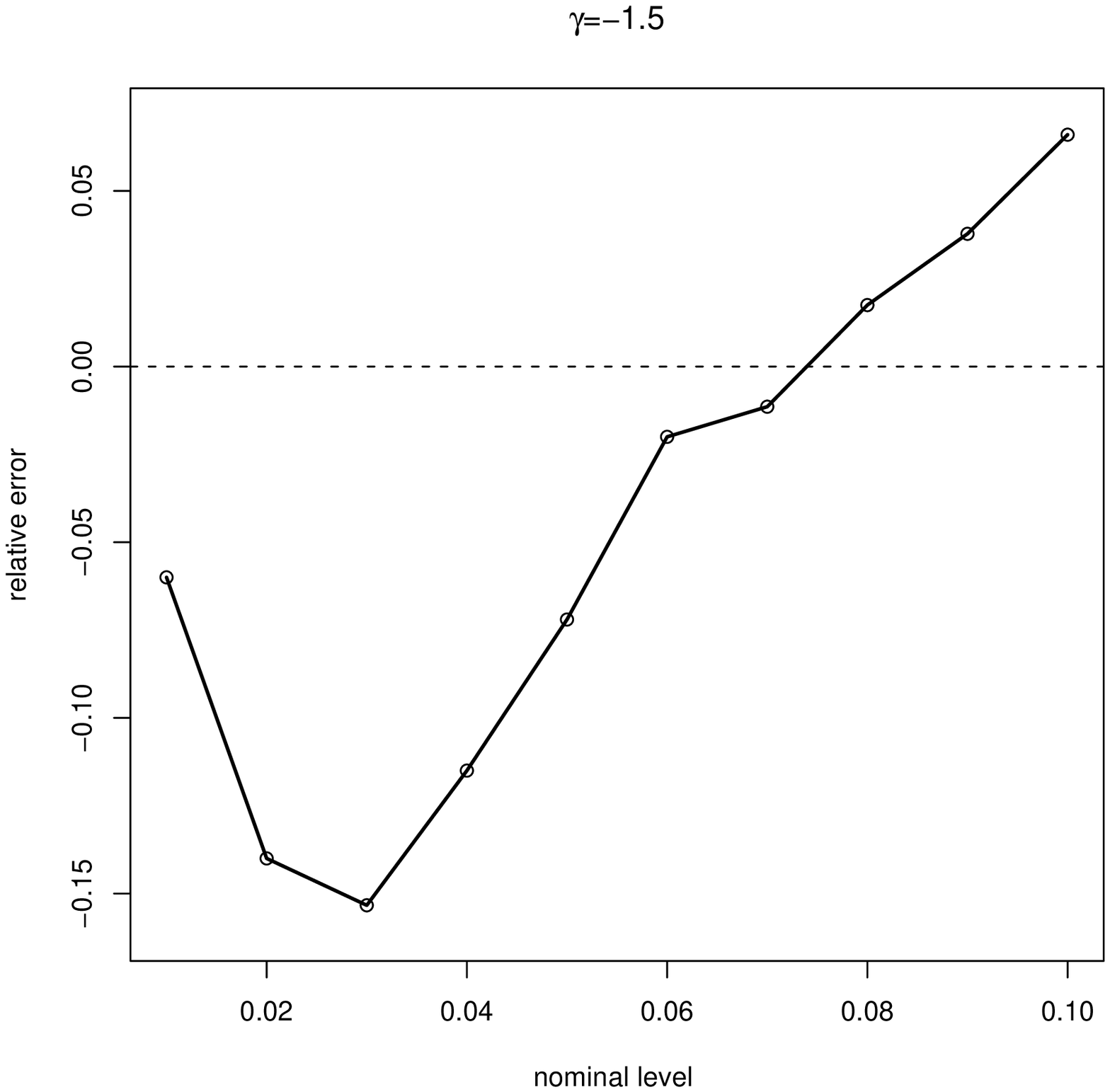}
\end{tabular}
\vspace{-2mm}} \caption{Relative errors of the robust tests
applied to the scale normal model $\mathcal{N}(0,1)$, when
$\overline{\sigma}=1.9$, 98 data are generated from model and 2
outliers x=10 are added.}\label{Gr8}
\end{figure}
}
\end{center}

\begin{center}
{\scriptsize
\begin{figure}[tbp]
{\scriptsize
\begin{tabular}{cc}
\includegraphics[width=6.5cm,height=6cm]{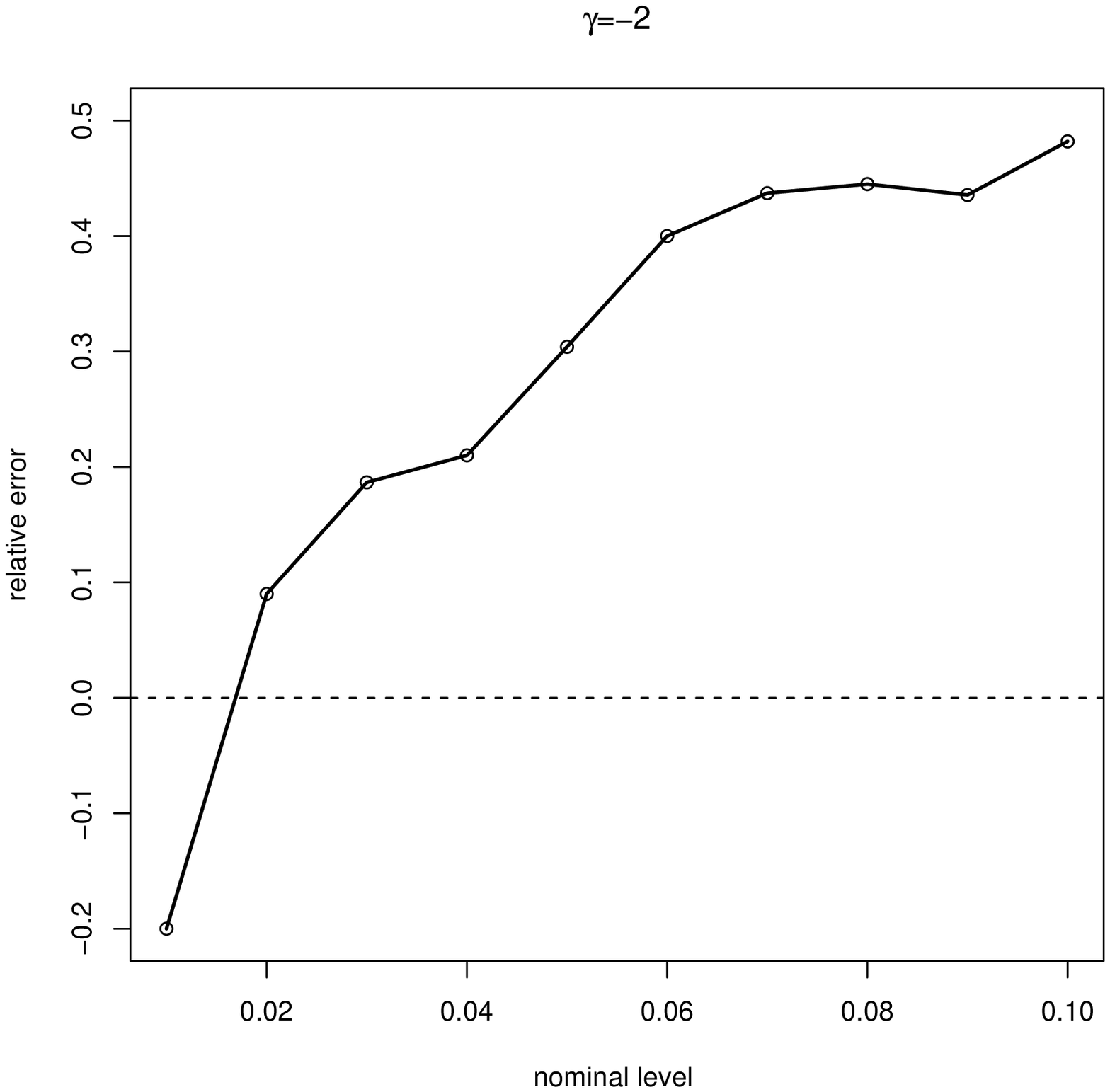} &
\includegraphics[width=6.5cm,height=6cm]{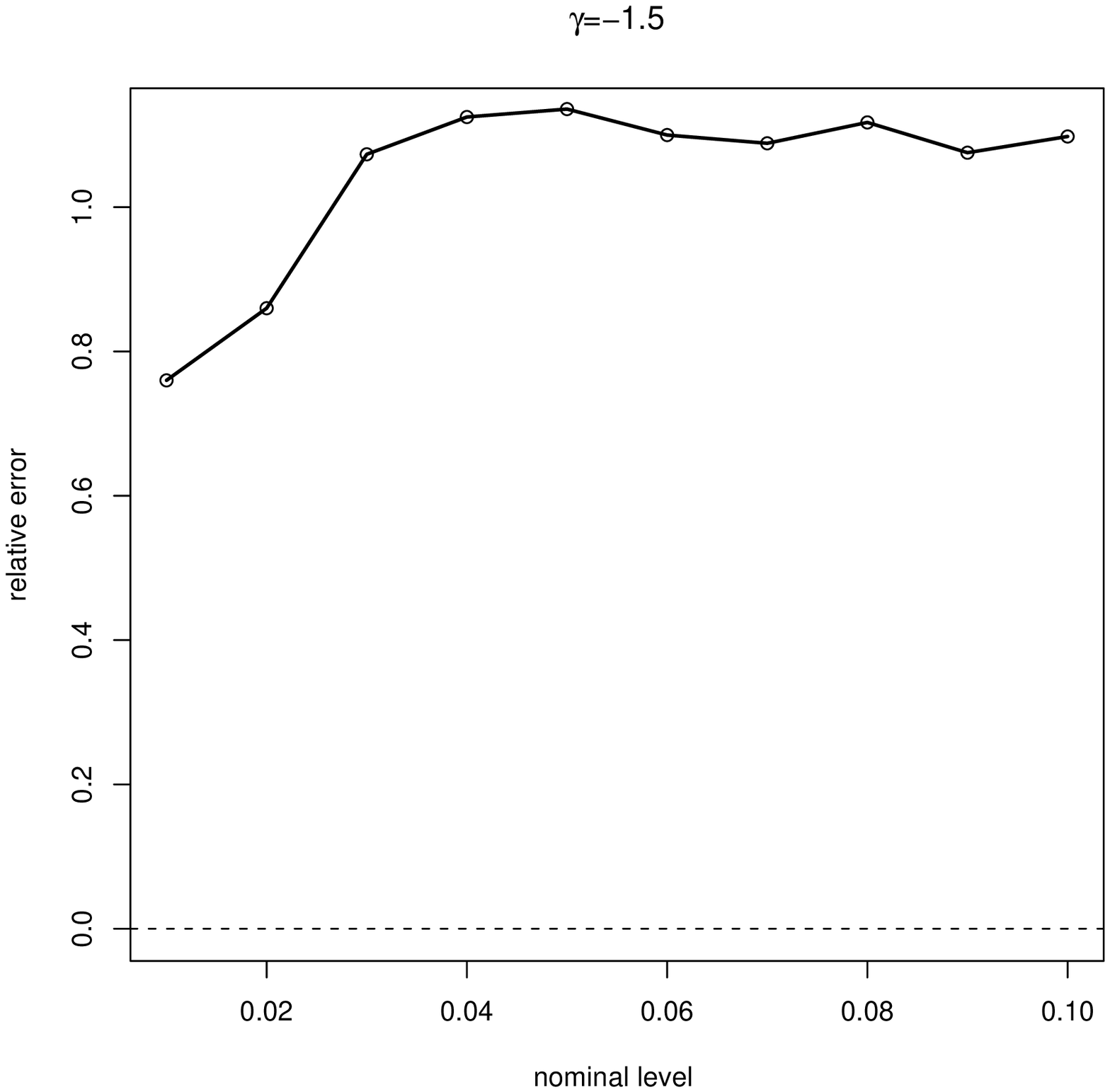}
\end{tabular}
\vspace{-2mm}} \caption{Relative errors of the robust tests
applied to the scale normal model $\mathcal{N}(0,1)$, when
$\overline{\sigma}=1.9$, 96 data are generated from model and 4
outliers x=10 are added.}\label{Gr9}
\end{figure}
}
\end{center}

\newpage

\newpage

\vspace{10cm}

\section{An adaptive choice of the tuning
parameter}\label{secchoice}

At the present stage we only present some heuristic and defer the
formal treatment of this proposal, which lays beyond the scope of
the present work.

According to the model and the parameter to be estimated, the
choice of $\gamma$ should be considered with respect to the
expression (\ref{term}) which has to be bounded. We refer to the
examples given in subsection \ref{study} for some scale and
location model.

Given a set of observations $X_{1},\dots,X_{n}$ an adaptive choice
for $\alpha$ would aim at reducing the estimated maximal bias
caused by an extraneous data.  Define
$\widehat{\theta}_{n}(\alpha,\gamma)$ the D$\phi$E of $\theta_{0}$
on the entire set of observation. For $1\leq i\leq n$, let
$\widehat{\theta}_{n-1}^{i}(\alpha,\gamma)$ be the D$\phi$E of
$\theta_{0}$
built on the leave one out data set $X_{1},X_{2},\dots,X_{i-1},\newline X_{i+1}%
,\dots,X_{n}$. Define
\[
B_{n}(\alpha,\gamma):=\max_{i}|\widehat{\theta}_{n}(\alpha,\gamma
)-\widehat{\theta}_{n-1}^{i}(\alpha,\gamma)|
\]
which measures the maximal bias caused by a single outlier and
\[
\alpha^{\ast}(\gamma):=\arg\inf_{\alpha}B_{n}(\alpha,\gamma).
\]

\section{Proofs}

\textit{Proof of Proposition \ref{t1b}}

For fixed $\alpha$, $\widehat{\theta}_{n}(\alpha)$ are
M-estimators. In accordance with the theory regarding the
M-estimators (see for example van der Vaart \cite{va}), the so
called $\psi$-function corresponding to
$\widehat{\theta}_{n}(\alpha)$ is
\begin{equation*}
\psi_{\alpha}(x,\theta)=m^{\prime}(\theta,\alpha,x)
\end{equation*}
and the influence function of $T_{\alpha}$ is
\begin{equation}
\label{e5b}
\mathrm{IF}(x;T_{\alpha},P_{\theta_{0}})=\left[M(\psi_{\alpha},P_{\theta_{0}})\right]^{-1}\psi_{\alpha}\left(x,T_{\alpha}(P_{\theta_{0}})\right)
\end{equation}
where
\begin{equation*}
M\left(\psi_{\alpha},P_{\theta_{0}}\right)=-\int\frac{\partial}{\partial
\theta}\left[\psi_{\alpha}(y,\theta)\right]_{\theta_{0}}dP_{\theta_{0}}(y)=-\int
m^{\prime \prime}(\theta_{0},\alpha,y)dP_{\theta_{0}}(y).
\end{equation*}
Using the Fisher consistency of the functional $T_{\alpha}$,
\begin{eqnarray*}
\psi_{\alpha}\left(x,T_{\alpha}(P_{\theta_{0}})\right) &=&
\psi_{\alpha}(x,\theta_{0})\\ &=&
-\int\varphi^{\prime\prime}\left(\frac{p_{\alpha}}{p_{\theta_{0}}}\right)\frac{p_{\alpha}}{p_{\theta_{0}}^{2}}\dot{p}_{\theta_{0}}dP_{\alpha}+\varphi^{\prime\prime}\left(\frac{p_{\alpha}}{p_{\theta_{0}}}(x)\right)\frac{p_{\alpha}^{2}(x)}{p_{\theta_{0}}^{3}(x)}\dot{p}_{\theta_{0}}(x)
\end{eqnarray*}
which substituted in (\ref{e5b}) leads to the announced
result.\hfill $\square$

\textit{Proof of Proposition \ref{t2b}}

Let $\varepsilon>0$ and $\widetilde{P_{\theta_{0}}}_{\varepsilon
x}=(1-\varepsilon)P_{\theta_{0}}+\varepsilon \delta_{x}$ be the
contaminated model. Then
\begin{eqnarray*}
 U_{\alpha}\left(\widetilde{P_{\theta_{0}}}_{\varepsilon
 x}\right) &=&\int m\left(T_{\alpha}\left(\widetilde{P_{\theta_{0}}}_{\varepsilon
 x}\right),\alpha,y\right)d\widetilde{P_{\theta_{0}}}_{\varepsilon
x}(y)\\ &=& (1-\varepsilon)\int
m\left(T_{\alpha}\left(\widetilde{P_{\theta_{0}}}_{\varepsilon
 x}\right),\alpha,y\right)dP_{\theta_{0}}(y)+\varepsilon m\left(T_{\alpha}\left(\widetilde{P_{\theta_{0}}}_{\varepsilon
 x}\right),\alpha,x\right)
 \end{eqnarray*}
and derivation yields
 \begin{eqnarray*} && \mathrm{IF}(x;U_{\alpha},P_{\theta_{0}})=\frac{\partial}{\partial\varepsilon}\left[ U_{\alpha}\left(\widetilde{P_{\theta_{0}}}_{\varepsilon
 x}\right)\right]_{\varepsilon=0}=\\
  &=&-\int m(\theta_{0},\alpha,y)dP_{\theta_{0}}(y)+\mathrm{IF}(x;T_{\alpha},P_{\theta_{0}})^{t}\int
 m^{\prime}(\theta_{0},\alpha,y)dP_{\theta_{0}}(y)+m(\theta_{0},\alpha,x)\\ &=&
-\phi(\alpha,\theta_{0})+m(\theta_{0},\alpha,x).
 \end{eqnarray*}
\hfill$\square$

\textit{Proof of Proposition \ref{t3b}}

For notational clearness, define $T:\Theta\times \mathcal{M}\to
\Theta$,
\begin{equation*}T(\alpha,P):=T_{\alpha}(P).\end{equation*}

For each $\alpha\in\Theta$, the definition of $T(\alpha,P)$ leads
to
\begin{equation*}
\int m^{\prime}\left(T(\alpha,P),\alpha,y\right)dP(y)=0.
\end{equation*}

By the very definition of $V(P)$ and $T(V(P),P)$, they both obey
\begin{equation}
\label{e12b} \left\{ \begin{array}{ll}
                     \int
m^{\prime}\left(T(V(P),P),V(P),y\right)dP(y)=0\\
                    \int\frac{\partial}{\partial
                    \alpha}\left[m\left(T(\alpha,P),\alpha,y\right)\right]_{V(P)}dP(y)=0
                    \end{array}.
                    \right.
\end{equation}
Denoting $n(\theta,\alpha,y)=\frac{\partial}{\partial
 \alpha}m(\theta,\alpha,y)$
\begin{eqnarray*}
n(\theta,\alpha,y) &=&
\int\varphi^{\prime\prime}\left(\frac{p_{\alpha}}{p_{\theta}}\right)\frac{\dot{p}_{\alpha}}{p_{\theta}}dP_{\alpha}+\int\varphi^{\prime}\left(\frac{p_{\alpha}}{p_{\theta}}\right)\frac{\dot{p}_{\alpha}}{p_{\alpha}}dP_{\alpha}-\\&&
-\left\{\varphi^{\prime
\prime}\left(\frac{p_{\alpha}}{p_{\theta}}(y)\right)\frac{p_{\alpha}}{p_{\theta}^{2}}(y)\dot{p}_{\alpha}(y)
+\varphi^{\prime}\left(\frac{p_{\alpha}}{p_{\theta}}(y)\right)\frac{\dot{p}_{\alpha}(y)}{p_{\theta}(y)}-\varphi^{\prime}\left(\frac{p_{\alpha}}{p_{\theta}}(y)\right)\frac{\dot{p}_{\alpha}(y)}{p_{\theta}(y)}\right\}\\
 &=& \int\left\{\varphi^{\prime
 \prime}\left(\frac{p_{\alpha}}{p_{\theta}}\right)\frac{1}{p_{\theta}}+\varphi^{\prime}\left(\frac{p_{\alpha}}{p_{\theta}}\right)\frac{1}{p_{\alpha}}\right\}\dot{p}_{\alpha}dP_{\alpha}-\varphi^{\prime
\prime}\left(\frac{p_{\alpha}}{p_{\theta}}(y)\right)\frac{p_{\alpha}}{p_{\theta}^{2}}(y)\dot{p}_{\alpha}(y).
 \end{eqnarray*}
 From (\ref{e12b})
 \[\left\{ \begin{array}{ll}
                     \int
m^{\prime}\left(T(V(P),P),V(P),y\right)dP(y)=0\\
\!\frac{\partial}{\partial\alpha}\left[T(\alpha,P)\right]_{V(P)}\!\!\int\!
m^{\prime}\left(T(V(P),P),V(P),y\right)dP(y)\!+\!\!\int
\!n\left(T(V(P),P),V(P),y\right)\!dP(y)\!=\!0
                    \end{array}
                    \right.\]
and consequently
 \[\left\{ \begin{array}{ll}
                     \int
m^{\prime}\left(T(V(P),P),V(P),y\right)dP(y)=0\\ \int
n\left(T(V(P),P),V(P),y\right)dP(y)=0
                    \end{array}.
                    \right.\]

For the contaminated model
 \begin{equation*}
 \int
 n\left(T\left(V\left(\widetilde{P_{\theta_{0}}}_{\varepsilon x}\right),\widetilde{P_{\theta_{0}}}_{\varepsilon x}\right),V\left(\widetilde{P_{\theta_{0}}}_{\varepsilon x}\right),y\right)d\widetilde{P_{\theta_{0}}}_{\varepsilon x}(y)=0
 \end{equation*}
 and so
\[ \begin{array}{ll}
 (1-\varepsilon) \int
 n\left(T\left(V\left(\widetilde{P_{\theta_{0}}}_{\varepsilon x}\right),\widetilde{P_{\theta_{0}}}_{\varepsilon x}\right),V\left(\widetilde{P_{\theta_{0}}}_{\varepsilon x}\right),y\right)dP_{\theta_{0}}(y)+\\+\varepsilon
 n\left(T\left(V\left(\widetilde{P_{\theta_{0}}}_{\varepsilon x}\right),\widetilde{P_{\theta_{0}}}_{\varepsilon x}\right),V\left(\widetilde{P_{\theta_{0}}}_{\varepsilon
 x}\right),x\right)=0.
\end{array}
\]
Now derivation yields \begin{equation}\label{e13b}
\begin{array}{cc} -\int
n(\theta_{0},\theta_{0},y)dP_{\theta_{0}}(y)+\int\frac{\partial}{\partial
\theta
}n\left(T(\theta_{0},P_{\theta_{0}}),\theta_{0},y\right)dP_{\theta_{0}}(y)\left\{\frac{\partial}{\partial
\alpha
}T(\theta_{0},P_{\theta_{0}})\mathrm{IF}(x;V,P_{\theta_{0}})+\right.\\
\left. +\mathrm{IF}(x;T_{\theta_{0}},P_{\theta_{0}})\right\}
+\int\frac{\partial}{\partial
\alpha}n(\theta_{0},\theta_{0},y)dP_{\theta_{0}}(y)\mathrm{IF}(x;V,P_{\theta_{0}})+n(\theta_{0},\theta_{0},x)=0.
\end{array}
\end{equation}

Since
\begin{eqnarray*}
\frac{\partial}{\partial \theta}n(\theta,\alpha,y)&=&
\int\left\{-2\varphi^{\prime
\prime}\left(\frac{p_{\alpha}}{p_{\theta}}\right)\frac{\dot{p}_{\theta}}{p_{\theta}^{2}}-\varphi^{\prime
\prime
\prime}\left(\frac{p_{\alpha}}{p_{\theta}}\right)\frac{p_{\alpha}}{p_{\theta}^{3}}\dot{p}_{\theta}\right\}\dot{p}_{\alpha}^{t}dP_{\alpha}+\\
&&+\left\{\varphi^{\prime \prime
\prime}\left(\frac{p_{\alpha}}{p_{\theta}}(y)\right)+2\varphi^{\prime
\prime}\left(\frac{p_{\alpha}}{p_{\theta}}(y)\right)\frac{p_{\theta}(y)}{p_{\alpha}(y)}\right\}\frac{p_{\alpha}^{2}(y)}{p_{\theta}^{4}(y)}\dot{p}_{\theta}(y)\dot{p}_{\alpha}(y)^{t}
\end{eqnarray*}
and particularly
\begin{equation*}
\frac{\partial}{\partial
\theta}n(\theta_{0},\theta_{0},y)=\left\{2\varphi^{\prime\prime}(1)+\varphi^{\prime
\prime\prime}(1)\right\}\left\{\frac{\dot{p}_{\theta_{0}}(y)\dot{p}_{\theta_{0}}(y)^{t}}{p_{\theta_{0}}^{2}(y)}-\int\frac{\dot{p}_{\theta_{0}}\dot{p}_{\theta_{0}}^{t}}{p_{\theta_{0}}^{2}}dP_{\theta_{0}}\right\},
\end{equation*}
deduce that
\begin{equation}
\label{e14b} \int\frac{\partial}{\partial
\theta}n(\theta_{0},\theta_{0},y)dP_{\theta_{0}}(y)=0.
\end{equation}

On the other hand
\begin{eqnarray*}
\frac{\partial}{\partial
\alpha}n(\theta,\alpha,y)&=&\int\left\{\varphi^{\prime\prime\prime}\left(\frac{p_{\alpha}}{p_{\theta}}\right)\frac{\dot{p}_{\alpha}}{p_{\theta}^{2}}-\varphi^{\prime}\left(\frac{p_{\alpha}}{p_{\theta}}\right)\frac{\dot{p}_{\alpha}}{p_{\alpha}^{2}}+\varphi^{\prime\prime}\left(\frac{p_{\alpha}}{p_{\theta}}\right)\frac{\dot{p}_{\alpha}}{p_{\alpha}p_{\theta}}\right\}\dot{p}_{\alpha}^{t}dP_{\alpha}+\\
&&+\int\left\{\varphi^{\prime\prime}\left(\frac{p_{\alpha}}{p_{\theta}}\right)\frac{1}{p_{\theta}}+\varphi^{\prime}\left(\frac{p_{\alpha}}{p_{\theta}}\right)\frac{1}{p_{\alpha}}\right\}\left(\ddot{p}_{\alpha}+\frac{\dot{p}_{\alpha}\dot{p}_{\alpha}^{t}}{p_{\alpha}}\right)dP_{\alpha}-\\
&&-\varphi^{\prime\prime\prime}\left(\frac{p_{\alpha}}{p_{\theta}}(y)\right)\frac{p_{\alpha}(y)}{p_{\theta}^{3}(y)}\dot{p}_{\alpha}(y)\dot{p}_{\alpha}(y)^{t}-\\
&&-\varphi^{\prime\prime}\left(\frac{p_{\alpha}}{p_{\theta}}(y)\right)\frac{\dot{p}_{\alpha}(y)\dot{p}_{\alpha}(y)^{t}}{p_{\theta}^{2}(y)}-\varphi^{\prime\prime}\left(\frac{p_{\alpha}}{p_{\theta}}(y)\right)\frac{p_{\alpha}(y)}{p_{\theta}^{2}(y)}\ddot{p}_{\alpha}(y)
\end{eqnarray*}
and particularly
\begin{eqnarray*}
\frac{\partial}{\partial
\alpha}n(\theta_{0},\theta_{0},y)&=&\int\left\{\varphi^{\prime\prime\prime}(1)+2\varphi^{\prime\prime}(1)\right\}\frac{\dot{p}_{\theta_{0}}\dot{p}_{\theta_{0}}^{t}}{p_{\theta_{0}}^{2}}dP_{\theta_{0}}+\int
\varphi^{\prime\prime}(1)
\frac{\ddot{p}_{\theta_{0}}}{p_{\theta_{0}}}dP_{\theta_{0}}\\
&&-\{\varphi^{\prime\prime\prime}(1)+\varphi^{\prime\prime}(1)\}\frac{\dot{p}_{\theta_{0}}(y)\dot{p}_{\theta_{0}}(y)^{t}}{p_{\theta_{0}}^{2}(y)}-\varphi^{\prime\prime}(1)\frac{\ddot{p}_{\theta_{0}}(y)}{p_{\theta_{0}}(y)}.
\end{eqnarray*}

As a consequence
\begin{equation}
\label{e15b} \int\frac{\partial}{\partial
\alpha}n(\theta_{0},\theta_{0},y)dP_{\theta_{0}}(y)=\varphi^{\prime\prime}(1)\int\frac{\dot{p}_{\theta_{0}}\dot{p}_{\theta_{0}}^{t}}{p_{\theta_{0}}^{2}}dP_{\theta_{0}}=\varphi^{\prime\prime}(1)I_{\theta_{0}}.
\end{equation}

Also
\begin{equation}
\label{e16b}
n(\theta_{0},\theta_{0},x)=-\varphi^{\prime\prime}(1)\frac{\dot{p}_{\theta_{0}}(x)}{p_{\theta_{0}}(x)}.
\end{equation}

Using the Fisher consistency of the functional $T_{\alpha}$ and
substituting (\ref{e14b}), (\ref{e15b}) and (\ref{e16b}) in
(\ref{e13b}) it holds
\begin{equation}
\label{e17b}
\mathrm{IF}(x;V,P_{\theta_{0}})=I_{\theta_{0}}^{-1}\frac{\dot{p}_{\theta_{0}}(x)}{p_{\theta_{0}}(x)}
\end{equation}
and this completes the proof.\hfill $\square$

\textit{Proof of Proposition \ref{psc}}

By replacing $\dot{p}_{\theta}$,
\begin{eqnarray}
\!\!\!\!\left(\frac{p_{\alpha}(x)}{p_{\theta_{0}}(x)}\right)^{\gamma}\dot{p}_{\theta_{0}}(x)
&=&
-\frac{\theta_{0}^{\gamma-2}}{\alpha^{\gamma}}\left\{\frac{p(\alpha^{-1}x)^{\gamma}}{p(\theta_{0}^{-1}x)^{\gamma-1}}+\frac{x}{\theta_{0}}\left(\frac{p(\alpha^{-1}x)}{p(\theta_{0}^{-1}x)}\right)^{\gamma}\dot{p}(\theta_{0}^{-1}x)\right\}\nonumber\\
&=&
-\frac{\theta_{0}^{\gamma-2}}{\alpha^{\gamma}}\left\{\frac{p(\alpha^{-1}x)^{\gamma}}{p(\theta_{0}^{-1}x)^{\gamma-1}}-\theta_{0}\frac{p(\alpha^{-1}x)^{\gamma}}{p(\theta_{0}^{-1}x)^{\gamma-1}}\frac{\partial}{\partial\theta}\log
p(\theta_{0}^{-1}x) \right\}\label{rht1}
\end{eqnarray}
and similarly
\begin{eqnarray*}
\!\!\!\left(\frac{p_{\alpha}(x)}{p_{\theta_{0}}(x)}\right)^{\gamma}\frac{\dot{p}_{\theta_{0}}(x)}{p_{\theta_{0}}(x)}
&=&-\frac{\theta_{0}^{\gamma-1}}{\alpha^{\gamma}}\left\{\left(\frac{p(\alpha^{-1}x)}{p(\theta_{0}^{-1}x)}\right)^{\gamma}+\frac{x}{\theta_{0}}\left(\frac{p(\alpha^{-1}x)}{p(\theta_{0}^{-1}x)}\right)^{\gamma}\frac{\dot{p}(\theta_{0}^{-1}x)}{p(\theta_{0}^{-1}x)}\right\}\\
&=&-\frac{\theta_{0}^{\gamma-1}}{\alpha^{\gamma}}\left\{\left(\frac{p(\alpha^{-1}x)}{p(\theta_{0}^{-1}x)}\right)^{\gamma}-\theta_{0}\left(\frac{p(\alpha^{-1}x)}{p(\theta_{0}^{-1}x)}\right)^{\gamma}\frac{\partial}{\partial\theta}\log
p(\theta_{0}^{-1}x) \right\}.
\end{eqnarray*}
The condition (A.1) together with one of the conditions (A.2) or
(A.3) (depending on the choice of $\gamma$) entails that the
function (\ref{rht1}) is integrable. On the other hand (A.2) or
(A.3) together with (A.4) assure that the function in the above
display is bounded. Then
$\mathrm{IF}(x;T_{\alpha},P_{\theta_{0}})$ as it is expressed by
(\ref{b1}) is a bounded function.\hfill $\square$

\textit{Proof of Proposition \ref{loccondd}}

It holds
\begin{equation*}
\left(\frac{p_{\alpha}(x)}{p_{\theta_{0}}(x)}\right)^{\gamma}\dot{p}_{\theta_{0}}(x)=\frac{p(x-\alpha)^{\gamma}}{p(x-\theta_{0})^{\gamma-1}}\frac{\partial}{\partial\theta}\log
p(x-\theta_{0})
\end{equation*}
and
\begin{equation*}
\left(\frac{p_{\alpha}(x)}{p_{\theta_{0}}(x)}\right)^{\gamma}\frac{\dot{p}_{\theta_{0}}(x)}{p_{\theta_{0}}(x)}
=\left(\frac{p(x-\alpha)}{p(x-\theta_{0})}\right)^{\gamma}\frac{\partial}{\partial\theta}\log
p(x-\theta_{0}).
\end{equation*}

Then the condition (\ref{loccond}) allows to conclude that
$\mathrm{IF}(x;T_{\alpha},P_{\theta_{0}})$, as it is expressed by
(\ref{b1}), is bounded.\hfill $\square$

\textit{Proof of Proposition \ref{proppas0}}

First prove that
$\widehat{\theta}_{n}(\alpha)-T_{\alpha}(P_{n,\varepsilon,x}^{P})=o_{P}(1).$

It holds $\int m^{\prime}(\theta_{0}, \alpha)dP_{\theta_{0}}=0$
and
\begin{equation}\label{SS}
\int m^{\prime\prime}(\theta_{0}, \alpha)dP_{\theta_{0}}=-\int
\varphi^{\prime\prime}\left(\frac{p_{\alpha}}{p_{\theta_{0}}}\right)
\frac{p_{\alpha}^{2}}{p_{\theta_{0}}^{3}}\dot{p}_{\theta_{0}}\dot{p}_{\theta_{0}}^{t}d\lambda=-S.
\end{equation}
The matrix $S$ is symmetric and positive since
$\varphi^{\prime\prime}$ is positive by the convexity of
$\varphi$. Using (C.3) in connection with the Lindeberg-Feller
theorem for triangular arrays we have $\sqrt{n}\int
m^{\prime}(\theta_{0}, \alpha)dP_{n}=O_{P}(1)$. Using (C.4) in
connection with the Lindeberg-Feller theorem for triangular arrays
yields $\int m^{\prime\prime}(\theta_{0},
\alpha)dP_{n}+S=o_{P}(1)$.

Now, for any $\theta=\theta_{0}+un^{-1/3}$ with $\|u\|\leq 1$, a
Taylor expansion of $\int m(\theta, \alpha)dP_{n}$ around
$\theta_{0}$ under (C.1) yields
\begin{eqnarray*}
&&n\int m(\theta, \alpha)dP_{n}-n \int m(\theta_{0},
\alpha)dP_{n}=\\&&=n^{2/3}u^{t} \int m^{\prime}(\theta_{0},
\alpha)dP_{n}+2^{-1}n^{1/3}u^{t}\int m^{\prime\prime}(\theta_{0},
\alpha)dP_{n}u+O_{P}(1)
\end{eqnarray*}
uniformly on $u$ with $\|u\|\leq 1$. Hence
\begin{equation*}
n\int m(\theta, \alpha)dP_{n}-n \int m(\theta_{0},
\alpha)dP_{n}=O_{P}(n^{1/6})-2^{-1}n^{1/3}u^{t}Su+O_{P}(1)
\end{equation*}
uniformly on $u$ with $\|u\|\leq 1$. Hence uniformly on $u$ with
$\|u\|=1$,
\begin{equation}\label{inegg}
n\int m(\theta, \alpha)dP_{n}-n \int m(\theta_{0},
\alpha)dP_{n}\leq O_{P}(n^{1/6})-2^{-1}cn^{1/3}+O_{P}(1)
\end{equation}
where $c$ is the smallest eigenvalue of the matrix $S$. Note that
$c$ is positive since $S$ is positive definite. In view of
(\ref{inegg}), by the continuity of $\theta\to \int m(\theta,
\alpha)dP_{n}$, it holds that as $n\to\infty$, with probability
one, $\theta\to \int m(\theta, \alpha)dP_{n}$ attains its maximum
at some point $\widehat{\theta}_{n}(\alpha)$ in the interior of
the ball $\{\theta: \|\theta-\theta_{0}\|\leq n^{-1/3}\}$, and
therefore
\begin{equation}\label{CAm}
\widehat{\theta}_{n}(\alpha)-\theta_{0}=o_{P}(1).
\end{equation}

On the other hand,
\begin{equation*}
T_{\alpha}(P_{n,\varepsilon,x}^{P})=\theta_{0}+\frac{\varepsilon}{\sqrt{n}}\mathrm{IF}(x;T_{\alpha},P_{\theta_{0}})+\frac{\Delta}{\sqrt{n}}{\bf{1}}+\rho\left(\frac{\varepsilon}{\sqrt{n}},
\frac{\Delta}{\sqrt{n}}\right)\frac{\sqrt{\varepsilon^{2}+\Delta^{2}}}{n}
\end{equation*}
where ${\bf{1}}:=(1,\dots,1)^{t}$ above coincides with
$\frac{\partial}{\partial
\widetilde{\Delta}}[T_{\alpha}(P_{\theta_{0}+\widetilde{\Delta}})]_{\widetilde{\Delta}=0}$
by the Fisher consistency of the functional $T_{\alpha}$ and the
function $\rho$ satisfies $\lim_{n\to
\infty}\rho\left(\frac{\varepsilon}{\sqrt{n}},
\frac{\Delta}{\sqrt{n}}\right)=0$. Then
$T_{\alpha}(P_{n,\varepsilon,x}^{P})-\theta_{0}$ converges to zero
in probability as $n\to\infty$. Combining this with (\ref{CAm}) we
obtain that
$\widehat{\theta}_{n}(\alpha)-T_{\alpha}(P_{n,\varepsilon,x}^{P})$
converges to zero in probability.

In the following, we prove that
$\sqrt{n}(\widehat{\theta}_{n}(\alpha)-T_{\alpha}(P_{n,\varepsilon,x}^{P}))=O_{P}(1)$.

By Taylor expansion, there exists $\overline{\theta}_{n}$ inside
the segment that links $T_{\alpha}(P_{n,\varepsilon,x}^{P})$ and
$\widehat{\theta}_{n}(\alpha)$ such that
\begin{eqnarray}
0&=& \int
m^{\prime}(\widehat{\theta}_{n}(\alpha),\alpha)dP_{n}\nonumber\\
&=& \int
m^{\prime}(T_{\alpha}(P_{n,\varepsilon,x}^{P}),\alpha)dP_{n}+\int
m^{\prime\prime}(T_{\alpha}(P_{n,\varepsilon,x}^{P}),\alpha)dP_{n}(\widehat{\theta}_{n}(\alpha)-T_{\alpha}(P_{n,\varepsilon,x}^{P}))+\nonumber\\
&&+\frac{1}{2}(\widehat{\theta}_{n}(\alpha)-T_{\alpha}(P_{n,\varepsilon,x}^{P}))^{t}\int
m^{\prime\prime\prime}(\overline{\theta}_{n},\alpha)dP_{n}(\widehat{\theta}_{n}(\alpha)-T_{\alpha}(P_{n,\varepsilon,x}^{P})).\label{as1}
\end{eqnarray}

By condition (C.1), using the sup-norm
\begin{equation*}
\|\int
m^{\prime\prime\prime}(\overline{\theta}_{n},\alpha)dP_{n}\|=\|\frac{1}{n}\sum_{k=1}^{n}m^{\prime\prime\prime}(\overline{\theta}_{n},\alpha)(X_{k})\|\leq\frac{1}{n}\sum_{k=1}^{n}H(X_{k}).
\end{equation*}
Applying the Lindeberg-Feller theorem for triangular arrays yields
$\int
m^{\prime\prime\prime}(\overline{\theta}_{n},\alpha)dP_{n}=O_{P}(1).$
Then the last term in (\ref{as1}) writes
$o_{P}(1)(\widehat{\theta}_{n}(\alpha)-T_{\alpha}(P_{n,\varepsilon,x}^{P}))$.

Under (C.1) and (C.4), by applying a Taylor expansion and
repeatedly the Lindeberg-Feller theorem for triangular arrays,
\begin{eqnarray*}
&&\int
m^{\prime\prime}(T_{\alpha}(P_{n,\varepsilon,x}^{P}),\alpha)dP_{n}
= \frac{1}{n}\sum
_{k=1}^{n}m^{\prime\prime}(T_{\alpha}(P_{n,\varepsilon,x}^{P}),\alpha,X_{k})\\
&=& \!\frac{1}{n}\sum
_{k=1}^{n}m^{\prime\prime}(\theta_{0},\alpha,X_{k})+\frac{\varepsilon}{n\sqrt{n}}\sum_{k=1}^{n}m^{\prime\prime\prime}(\theta_{0},\alpha,X_{k})\mathrm{IF}(x;T_{\alpha},P_{\theta_{0}})+\\
&&
+\frac{\Delta}{n\sqrt{n}}\sum_{k=1}^{n}m^{\prime\prime\prime}(\theta_{0},\alpha,X_{k}){\bf{1}}
+\rho\left(\frac{\varepsilon}{\sqrt{n}},
\frac{\Delta}{\sqrt{n}}\right)\frac{\sqrt{\varepsilon^{2}+\Delta^{2}}}{n}\\
&=& P_{\theta_{0}} m^{\prime\prime}(\theta_{0},\alpha)+o_{P}(1)
\end{eqnarray*}
where ${\bf{1}}:=(1,\dots,1)^{t}$ above coincides with
$\frac{\partial}{\partial
\widetilde{\Delta}}[T_{\alpha}(P_{\theta_{0}+\widetilde{\Delta}})]_{\widetilde{\Delta}=0}$
by the Fisher consistency of the functional $T_{\alpha}$ and the
function $\rho$ satisfies $\lim_{n\to
\infty}\rho\left(\frac{\varepsilon}{\sqrt{n}},
\frac{\Delta}{\sqrt{n}}\right)=0$.

Therefore (\ref{as1}) becomes
\begin{equation}\label{as2}
-\int
m^{\prime}(T_{\alpha}(P_{n,\varepsilon,x}^{P}),\alpha)dP_{n}=(\int
m^{\prime\prime}(\theta_{0},\alpha)dP_{\theta_{0}}+o_{P}(1))(\widehat{\theta}_{n}(\alpha)-T_{\alpha}(P_{n,\varepsilon,x}^{P})).
\end{equation}

We prove that $\sqrt{n}\int
m^{\prime}(T_{\alpha}(P_{n,\varepsilon,x}^{P}),\alpha)dP_{n}$ is
$O_{P}(1)$. By Taylor expansion,
\begin{eqnarray*}
\!\!&&\!\!\int
m^{\prime}(T_{\alpha}(P_{n,\varepsilon,x}^{P}),\alpha)dP_{n}\!=\frac{1}{n}\sum_{k=1}^{n}m^{\prime}(\theta_{0},\alpha,
X_{k}
)\!+\!\frac{\varepsilon}{n\sqrt{n}}\sum_{k=1}^{n}m^{\prime\prime}(\theta_{0},\alpha,
X_{k})\mathrm{IF}(x;T_{\alpha}, P_{\theta_{0}})+\\
&&+\frac{\Delta}{n\sqrt{n}}\sum_{k=1}^{n}m^{\prime\prime}(\theta_{0},\alpha,
X_{k}){\bf{1}} +\rho\left(\frac{\varepsilon}{\sqrt{n}},
\frac{\Delta}{\sqrt{n}}\right)\frac{\sqrt{\varepsilon^{2}+\Delta^{2}}}{n}
\end{eqnarray*}
and therefore
\begin{eqnarray*}
\!\!\!\!&&\sqrt{n}\int
m^{\prime}(T_{\alpha}(P_{n,\varepsilon,x}^{P}),\alpha)dP_{n}=\frac{1}{\sqrt{n}}\sum_{k=1}^{n}m^{\prime}(\theta_{0},\alpha,
X_{k}
)+\\\!\!\!\!&&+\frac{\varepsilon}{n}\sum_{k=1}^{n}m^{\prime\prime}(\theta_{0},\alpha,
X_{k})\mathrm{IF}(x;T_{\alpha},
P_{\theta_{0}})+\frac{\Delta}{n}\sum_{k=1}^{n}m^{\prime\prime}(\theta_{0},\alpha,
X_{k}){\bf{1}} +\rho\left(\frac{\varepsilon}{\sqrt{n}},
\frac{\Delta}{\sqrt{n}}\right)\sqrt{\frac{\varepsilon^{2}+\Delta^{2}}{n}}.
\end{eqnarray*}

Under (C.3) and (C.4), by applying the Lindeberg-Feller theorem
for triangular arrays it holds $\sqrt{n}\int
m^{\prime}(T_{\alpha}(P_{n,\varepsilon,x}^{P}),\alpha)dP_{n}=O_{P}(1)$.
Then from (\ref{as2})
\begin{equation*}\sqrt{n}(\widehat{\theta}_{n}(\alpha)-T_{\alpha}(P_{n,\varepsilon,x}^{P}))=O_{P}(1).\end{equation*}\hfill $\square$

\textit{Proof of Proposition \ref{proppas3}}

By Taylor expansion, there exists $\overline{\theta}_{n}$ inside
the segment that links $T_{\alpha}(P_{n,\varepsilon,x}^{P})$ and
$\widehat{\theta}_{n}(\alpha)$ such that
\begin{eqnarray*}
\!\!\!\!\!\!\!\!\!\!\!&&\widehat{\phi}_{n}(\alpha,\theta_{0}) =
\int m(\widehat{\theta}_{n}(\alpha),\alpha)dP_{n}\\
\!\!\!\!\!\!\!\!&&= \int
m(T_{\alpha}(P_{n,\varepsilon,x}^{P}),\alpha)dP_{n}+ \int
m^{\prime}(T_{\alpha}(P_{n,\varepsilon,x}^{P}),\alpha)^{t}dP_{n}(\widehat{\theta}_{n}(\alpha)-T_{\alpha}(P_{n,\varepsilon,x}^{P}))+\\
\!\!&&+
\frac{1}{2}(\widehat{\theta}_{n}(\alpha)-T_{\alpha}(P_{n,\varepsilon,x}^{P}))^{t}\int
m^{\prime\prime}(T_{\alpha}(P_{n,\varepsilon,x}^{P}),\alpha)dP_{n}(\widehat{\theta}_{n}(\alpha)-T_{\alpha}(P_{n,\varepsilon,x}^{P}))+\\
\!\!&&+\frac{1}{3!}\!\!\!\!\!\sum_{1\leq i,j,k\leq
d}\!\!\!\!\!\!(\widehat{\theta}_{n}(\alpha)\!-\!T_{\alpha}(P_{n,\varepsilon,x}^{P}))_{i}(\widehat{\theta}_{n}(\alpha)\!-\!T_{\alpha}(P_{n,\varepsilon,x}^{P}))_{j}(\widehat{\theta}_{n}(\alpha)\!-\!T_{\alpha}(P_{n,\varepsilon,x}^{P}))_{k}\int
\frac{\partial^{3}m(\overline{\theta}_{n},\alpha)}{\partial\theta_{i}\partial\theta_{j}\partial\theta_{k}}dP_{n}.
\end{eqnarray*}
Then
\begin{equation*}
\frac{\sqrt{n}(\widehat{\phi}_{n}(\alpha,\theta_{0})-U_{\alpha}(P_{n,\varepsilon,x}^{P}))}{[\int
\mathrm{IF}^{2}(y;U_{\alpha},P_{n,\varepsilon,x}^{P})dP_{n,\varepsilon,x}^{P}(y)]^{1/2}}
=
\frac{\sqrt{n}(\int
m(T_{\alpha}(P_{n,\varepsilon,x}^{P}),\alpha)dP_{n}-U_{\alpha}(P_{n,\varepsilon,x}^{P}))}{[\int
\mathrm{IF}^{2}(y;U_{\alpha},P_{n,\varepsilon,x}^{P})dP_{n,\varepsilon,x}^{P}(y)]^{1/2}}+
\end{equation*}
\vspace{-2mm}\begin{equation*} +\frac{(\int
m^{\prime}(T_{\alpha}(P_{n,\varepsilon,x}^{P}),\alpha)dP_{n})^{t}\sqrt{n}(\widehat{\theta}_{n}(\alpha)-T_{\alpha}(P_{n,\varepsilon,x}^{P}))}{[\int
\mathrm{IF}^{2}(y;U_{\alpha},P_{n,\varepsilon,x}^{P})dP_{n,\varepsilon,x}^{P}(y)]^{1/2}}+\end{equation*}
\vspace{-2mm}\begin{equation*}+\frac{\sqrt{n}(\widehat{\theta}_{n}(\alpha)-T_{\alpha}(P_{n,\varepsilon,x}^{P}))^{t}\int
m^{\prime\prime}(T_{\alpha}(P_{n,\varepsilon,x}^{P}),\alpha)dP_{n}(\widehat{\theta}_{n}(\alpha)-T_{\alpha}(P_{n,\varepsilon,x}^{P}))}{2[\int
\mathrm{IF}^{2}(y;U_{\alpha},P_{n,\varepsilon,x}^{P})dP_{n,\varepsilon,x}^{P}(y)]^{1/2}}+
\end{equation*}
\vspace{-2mm}\begin{equation}
\label{dezv}\frac{\sqrt{n}\sum_{1\leq i,j,k\leq
d}(\widehat{\theta}_{n}(\alpha)\!-\!T_{\alpha}(P_{n,\varepsilon,x}^{P}))_{i}(\widehat{\theta}_{n}(\alpha)\!-\!T_{\alpha}(P_{n,\varepsilon,x}^{P}))_{j}(\widehat{\theta}_{n}(\alpha)\!-\!T_{\alpha}(P_{n,\varepsilon,x}^{P}))_{k}\int
\frac{\partial^{3}m(\overline{\theta}_{n},\alpha)}{\partial\theta_{i}\partial\theta_{j}\partial\theta_{k}}
dP_{n}}{3![\int
\mathrm{IF}^{2}(y;U_{\alpha},P_{n,\varepsilon,x}^{P})dP_{n,\varepsilon,x}^{P}(y)]^{1/2}}.
\end{equation}

In the following we analyze each term in the above display. It
holds
\begin{equation*}
\int
m(T_{\alpha}(P_{n,\varepsilon,x}^{P}),\alpha)dP_{n}-U_{\alpha}(P_{n,\varepsilon,x}^{P})=\frac{1}{n}\sum_{k=1}^{n}\{m(T_{\alpha}(P_{n,\varepsilon,x}^{P}),\alpha,
X_{k})-P_{n,\varepsilon,x}^{P}m(T_{\alpha}(P_{n,\varepsilon,x}^{P}),\alpha)\}.
\end{equation*}
Apply the Lindeberg-Feller theorem for the triangular array
\begin{equation*}
Z_{n,k}:=m(T_{\alpha}(P_{n,\varepsilon,x}^{P}),\alpha,
X_{k})-P_{n,\varepsilon,x}^{P}m(T_{\alpha}(P_{n,\varepsilon,x}^{P}),\alpha).
\end{equation*}
For this, compute first $\mathrm{Var}(Z_{n,k})$. Observe that
\begin{eqnarray*}
\mathrm{Var}(Z_{n,k}) &=&\int
m^{2}(T_{\alpha}(P_{n,\varepsilon,x}^{P}),\alpha,y)dP_{n,\varepsilon,x}^{P}(y)-\left(\int
m(T_{\alpha}(P_{n,\varepsilon,x}^{P}),\alpha,y)dP_{n,\varepsilon,x}^{P}(y)\right)^{2}=\\
&=& \left(1-\frac{\varepsilon}{\sqrt{n}}\right)\int
m^{2}(T_{\alpha}(P_{n,\varepsilon,x}^{P}),\alpha,y)dP_{\theta_{n}}(y)+\frac{\varepsilon}{\sqrt{n}}m^{2}(T_{\alpha}(P_{n,\varepsilon,x}^{P}),\alpha,x)-\\
&&-\left\{\left(1-\frac{\varepsilon}{\sqrt{n}}\right)\int
m(T_{\alpha}(P_{n,\varepsilon,x}^{P}),\alpha,y)dP_{\theta_{n}}(y)+\frac{\varepsilon}{\sqrt{n}}m(T_{\alpha}(P_{n,\varepsilon,x}^{P}),\alpha,x)\right\}^{2}.
\end{eqnarray*}
By Taylor expansions
\begin{eqnarray*}
m(T_{\alpha}(P_{n,\varepsilon,x}^{P}),\alpha,y) &=&
m(\theta_{0},\alpha,y)+\frac{\varepsilon}{\sqrt{n}}m^{\prime}(\theta_{0},\alpha,y)^{t}\mathrm{IF}(x;T_{\alpha},P_{\theta_{0}})+\\&&+\frac{\Delta}{\sqrt{n}}m^{\prime}(\theta_{0},\alpha,y)^{t}{\bf{1}}
+\rho\left(\frac{\varepsilon}{\sqrt{n}},
\frac{\Delta}{\sqrt{n}}\right)\frac{\sqrt{\varepsilon^{2}+\Delta^{2}}}{n}\\
m^{2}(T_{\alpha}(P_{n,\varepsilon,x}^{P}),\alpha,y)&=&
m^{2}(\theta_{0},\alpha,y)+2\frac{\varepsilon}{\sqrt{n}}m(\theta_{0},\alpha,y)m^{\prime}(\theta_{0},\alpha,y)^{t}\mathrm{IF}(x;T_{\alpha},P_{\theta_{0}})+\\&&2\frac{\Delta}{\sqrt{n}}m(\theta_{0},\alpha,y)m^{\prime}(\theta_{0},\alpha,y)^{t}{\bf{1}}
+\rho\left(\frac{\varepsilon}{\sqrt{n}},
\frac{\Delta}{\sqrt{n}}\right)\frac{\sqrt{\varepsilon^{2}+\Delta^{2}}}{n}.
\end{eqnarray*}
Hence the conditions (C.2) and (C.3) assure that
$\mathrm{Var(}Z_{n,k})$ is finite.

We now prove the equality
\begin{equation}\label{varznk}
\mathrm{Var(}Z_{n,k})=\int
\mathrm{IF}^{2}(y;U_{\alpha},P_{n,\varepsilon,x}^{P})dP_{n,\varepsilon,x}^{P}(y).
\end{equation}

By definition,
\begin{equation*}
\mathrm{IF}(y;U_{\alpha},P_{n,\varepsilon,x}^{P})=\frac{\partial}{\partial
t}[U_{\alpha}(\widetilde{P_{n,\varepsilon,x}^{P}}_{ty})]_{t=0},
\end{equation*}
where
$\widetilde{P_{n,\varepsilon,x}^{P}}_{ty}=(1-t)P_{n,\varepsilon,x}^{P}+t\delta_{y}$.
Also
\begin{equation*}
U_{\alpha}(\widetilde{P_{n,\varepsilon,x}^{P}}_{ty})=(1-t)\int
m(T_{\alpha}(\widetilde{P_{n,\varepsilon,x}^{P}}_{ty}),\alpha,z)dP_{n,\varepsilon,x}^{P}(z)+tm(T_{\alpha}(\widetilde{P_{n,\varepsilon,x}^{P}}_{ty}),\alpha,y)
\end{equation*}
whence
\begin{equation*}
\mathrm{IF}(y;U_{\alpha},P_{n,\varepsilon,x}^{P})=-\int
m(T_{\alpha}(P_{n,\varepsilon,x}^{P}),\alpha,z)dP_{n,\varepsilon,x}^{P}(z)+\end{equation*}
\vspace{-5mm}\begin{equation*}+\int
m^{\prime}(T_{\alpha}(P_{n,\varepsilon,x}^{P}),\alpha,z)dP_{n,\varepsilon,x}^{P}(z)\mathrm{IF}(y;T_{\alpha},P_{n,\varepsilon,x}^{P})+m(T_{\alpha}(P_{n,\varepsilon,x}^{P}),\alpha,y).
\end{equation*}
By the definition of $T_{\alpha}(P_{n,\varepsilon,x}^{P})$, it
holds $$\int
m^{\prime}(T_{\alpha}(P_{n,\varepsilon,x}^{P}),\alpha,z)dP_{n,\varepsilon,x}^{P}(z)=0$$
and hence (\ref{varznk}) holds. Here we observe that
$\mathrm{IF}(y;T_{\alpha}, P_{n,\varepsilon,x}^{P})$ is finite for
any $y$, for any $n$ and any $\Delta$, since
$\mathrm{IF}(y;T_{\alpha}, P_{\theta_{0}})$ is.

Thus, by Lindeberg-Feller theorem for triangular arrays, the first
term in the expansion (\ref{dezv}) converges in distribution to a
variable $\mathcal{N}(0,1)$.

We have $\int
m^{\prime}(T_{\alpha}(P_{n,\varepsilon,x}^{P}),\alpha)dP_{n}=o_{P}(1)$
since $\sqrt{n}\int
m^{\prime}(T_{\alpha}(P_{n,\varepsilon,x}^{P}),\alpha)dP_{n}=O_{P}(1)$
(see the proof of Proposition \ref{proppas0}). Also, it holds
$\int
m^{\prime\prime}(T_{\alpha}(P_{n,\varepsilon,x}^{P}),\alpha)dP_{n}=O_{P}(1)$
and $\int
\frac{\partial^{3}m(\overline{\theta}_{n},\alpha)}{\partial\theta_{i}\partial\theta_{j}\partial\theta_{k}}
dP_{n}=O_{P}(1)$.

Consequently, using  Proposition \ref{proppas0} we obtain the
announced result.\hfill $\square$

\textit{Proof of Proposition \ref{p11}}

The level $\alpha_{0}$ is given by
\begin{eqnarray*}
\alpha_{0}&=& 2P_{\theta_{0}}(\widehat{\phi}_{n}\geq
k_{n}(\alpha_{0}))\\ &=&
2P_{\theta_{0}}\left(\frac{\sqrt{n}(\widehat{\phi}_{n}-U_{\alpha}(P_{\theta_{0}}))}{[\int
\mathrm{IF}^{2}(y;U_{\alpha},P_{\theta_{0}})dP_{\theta_{0}}(y)
]^{1/2}}\geq
\frac{\sqrt{n}(k_{n}(\alpha_{0})-U_{\alpha}(P_{\theta_{0}}))}{[\int
\mathrm{IF}^{2}(y;U_{\alpha},P_{\theta_{0}})dP_{\theta_{0}}(y)
]^{1/2}}\right).
\end{eqnarray*}

Using the asymptotic normality of $\widehat{\phi}_{n}$ in the case
of uncontaminated observations (see Broniatowski and Keziou
\cite{BrKe07}),
\begin{equation*}\frac{
\sqrt{n}(k_{n}(\alpha_{0})-U_{\alpha}(P_{\theta_{0}}))}{[\int
\mathrm{IF}^{2}(y;U_{\alpha},P_{\theta_{0}})dP_{\theta_{0}}(y)
]^{1/2}}=\Phi^{-1}\left(1-\frac{\alpha_{0}}{2}\right)+o(1).
\end{equation*}
Therefore
\begin{equation}\label{t2}
k_{n}(\alpha_{0})=U_{\alpha}(P_{\theta_{0}})+\frac{1}{\sqrt{n}}\Phi^{-1}\left(1-\frac{\alpha_{0}}{2}\right)[\int
\mathrm{IF}^{2}(y;U_{\alpha},P_{\theta_{0}})dP_{\theta_{0}}(y)
]^{1/2}+o\left(\frac{1}{\sqrt{n}}\right).
\end{equation}

Now we are interested in the value of the asymptotic power, when
the underlying distribution deviates slightly from the model.
Using (\ref{t2})
\begin{eqnarray*}
\!&&P_{n,\varepsilon,x}=2P_{n,\varepsilon,x}^{P}(\widehat{\phi}_{n}\geq
k_{n}(\alpha_{0}) )\\ \!\!&=&
2P_{n,\varepsilon,x}^{P}\left(\frac{\sqrt{n}(\widehat{\phi}_{n}-U_{\alpha}(P_{n,\varepsilon,x}^{P}))}{[\int
\mathrm{IF}^{2}(y;U_{\alpha},P_{n,\varepsilon,x}^{P})dP_{n,\varepsilon,x}^{P}(y)
]^{1/2}}\geq
\frac{\sqrt{n}(k_{n}(\alpha_{0})-U_{\alpha}(P_{n,\varepsilon,x}^{P}))}{[\int
\mathrm{IF}^{2}(y;U_{\alpha},P_{n,\varepsilon,x}^{P})dP_{n,\varepsilon,x}^{P}(y)]^{1/2}}\right)\\
\!\!&=&
2P_{n,\varepsilon,x}^{P}\left(\frac{\sqrt{n}(\widehat{\phi}_{n}-U_{\alpha}(P_{n,\varepsilon,x}^{P}))}{[\int
\mathrm{IF}^{2}(y;U_{\alpha},P_{n,\varepsilon,x}^{P})dP_{n,\varepsilon,x}^{P}(y)
]^{1/2}}\geq
-\frac{\sqrt{n}(U_{\alpha}(P_{n,\varepsilon,x}^{P})-U_{\alpha}(P_{\theta_{0}}))}{[\int
\mathrm{IF}^{2}(y;U_{\alpha},P_{n,\varepsilon,x}^{P})dP_{n,\varepsilon,x}^{P}(y)]^{1/2}}+\right.\\
\!\!&&+\left.
\Phi^{-1}\left(1-\frac{\alpha_{0}}{2}\right)\frac{[\int
\mathrm{IF}^{2}(y;U_{\alpha},P_{\theta_{0}})dP_{\theta_{0}}(y)]^{1/2}}{[\int
\mathrm{IF}^{2}(y;U_{\alpha},P_{n,\varepsilon,x}^{P})dP_{n,\varepsilon,x}^{P}(y)]^{1/2}}+o(1)\right).
\end{eqnarray*}

Expand $U_{\alpha}(P_{n,\varepsilon,x}^{P})$ around to
$U_{\alpha}(P_{\theta_{0}})$ to obtain
\begin{eqnarray*}
\sqrt{n}(U_{\alpha}(P_{n,\varepsilon,x}^{P})-U_{\alpha}(P_{\theta_{0}}))&=&\varepsilon
\mathrm{IF}(x;U_{\alpha},P_{\theta_{0}})+\Delta
\frac{\partial}{\partial
\widetilde{\Delta}}[U_{\alpha}(P_{\theta_{0}+\widetilde{\Delta}})]_{\widetilde{\Delta}=0}+\\&&+\rho\left(\frac{\varepsilon}{\sqrt{n}},\frac{\Delta}{\sqrt{n}}\right)\sqrt{\frac{\varepsilon^{2}+\Delta^{2}}{n}}.
\end{eqnarray*}

Using the asymptotic normality of the test statistic when the
observations are i.i.d. with $P_{n,\varepsilon,x}^{P}$ and taking
into account that
\begin{equation*}
\lim_{n\to\infty}[\int
\mathrm{IF}^{2}(y;U_{\alpha},P_{n,\varepsilon,x}^{P})dP_{n,\varepsilon,x}^{P}(y)]^{1/2}=[\int
\mathrm{IF}^{2}(y;U_{\alpha},P_{\theta_{0}})dP_{\theta_{0}}(y)]^{1/2}
\end{equation*}
it holds
\begin{eqnarray}
\lim_{n\to
\infty}P_{n,\varepsilon,x}&=&2-2\Phi\left(\Phi^{-1}\left(1-\frac{\alpha_{0}}{2}\right)-\Delta\frac{\frac{\partial}{\partial
\widetilde{\Delta}}[U_{\alpha}(P_{\theta_{0}+\widetilde{\Delta}})]_{\widetilde{\Delta}=0}}{[\int
\mathrm{IF}^{2}(y;U_{\alpha},P_{\theta_{0}})dP_{\theta_{0}}(y)]^{1/2}}-\right.\nonumber\\&&-\left.\varepsilon\frac{\mathrm{IF}(x;U_{\alpha},P_{\theta_{0}})}{[\int
\mathrm{IF}^{2}(y;U_{\alpha},P_{\theta_{0}})dP_{\theta_{0}}(y)]^{1/2}}\right).\label{t3}
\end{eqnarray}
A simple calculation shows that $\frac{\partial}{\partial
\widetilde{\Delta}}[U_{\alpha}(P_{\theta_{0}+\widetilde{\Delta}})]_{\widetilde{\Delta}=0}=\int
m(\theta_{0},\alpha,y)\frac{\dot{p}_{\theta_{0}}(y)}{p_{\theta_{0}}(y)}dP_{\theta_{0}}(y)$.
Hence (\ref{asp}) holds. \hfill $\square$

\end{document}